\documentclass[oneside,12pt]{amsart}
\usepackage{amsmath, amsfonts,amsthm,times,graphics}

 \makeatletter
\renewcommand*\subjclass[2][2010]{%
  \def\@subjclass{#2}%
  \@ifundefined{subjclassname@#1}{%
    \ClassWarning{\@classname}{Unknown edition (#1) of Mathematics
      Subject Classification; using '2010'.}%
  }{%
    \@xp\let\@xp\subjclassname\csname subjclassname@#1\endcsname
  }%
}

 \makeatother

\theoremstyle{definition}

 \makeatletter
\renewcommand*\subjclass[2][2010]{%
  \def\@subjclass{#2}%
  \@ifundefined{subjclassname@#1}{%
    \ClassWarning{\@classname}{Unknown edition (#1) of Mathematics
      Subject Classification; using '1991\'.}%
  }{%
    \@xp\let\@xp\subjclassname\csname subjclassname@#1\endcsname
  }%
}
 \makeatother

\begin{document}

\title[Lucas' theorem:  
its generalizations, extensions...]{Lucas' theorem: 
its generalizations, extensions and applications (1878--2014)}
\author{Romeo Me\v strovi\' c}

\address{Maritime Faculty, University of Montenegro, Dobrota 36,
 85330 Kotor, Montenegro} \email{romeo@ac.me}

{\renewcommand{\thefootnote}{}
\setcounter{footnote}{0}}

\maketitle

  \begin{abstract} 
In 1878 \'{E}. Lucas proved a remarkable result which provides a simple way to compute
the binomial coefficient ${n\choose m}$ modulo a prime $p$ 
in terms of the binomial coefficients of the base-$p$ digits of 
$n$ and $m$:  {\it If $p$ is a prime,  $n=n_0+n_1p+\cdots +n_sp^s$ and 
$m=m_0+m_1p+\cdots +m_sp^s$ are the $p$-adic expansions of  nonnegative 
integers $n$ and $m$, then
 \begin{equation*}
{n\choose m}\equiv \prod_{i=0}^{s}{n_i\choose m_i}\pmod{p}.
  \end{equation*}}

The above congruence, the so-called {\it Lucas' theorem} 
(or {\it Theorem of Lucas}), 
plays an important role in Number Theory and Combinatorics. 
In this article, consisting of six sections, 
we provide a historical survey of Lucas type congruences, generalizations
of Lucas' theorem modulo prime powers,  Lucas like
theorems for some  generalized binomial coefficients,   
and some their applications.

In Section \ref{sec1} we present  the fundamental congruences modulo a prime
including the famous Lucas' theorem. In Section \ref{sec2} we 
mention several known proofs and some  consequences of Lucas' theorem. 
In Section \ref{sec3} we present 
a number of extensions and variations of Lucas' theorem modulo prime powers.  
In Section \ref{sec4} we consider the notions of the Lucas property and the double  
Lucas property, where we also present numerous integer sequences 
satisfying one of these properties or a certain Lucas type congruence. 
In Section \ref{sec5} we collect several known Lucas type congruences 
for some generalized binomial coefficients. In particular, this concerns
the Fibonomial coefficients, the Lucas $u$-nomial coefficients, 
the Gaussian $q$-nomial coefficients and their generalizations.
Finally, some applications of Lucas' theorem in Number Theory  and
Combinatorics  are given in Section \ref{sec6}. 
  \end{abstract} 

\vfill\eject

\vspace*{2cm}
\begin{center}
{\Large CONTENTS}
\end{center}
\vspace{1cm}
{\bf \ref{sec1}  Introduction} \hfill {\bf 3}\bigskip

\noindent{\bf \ref{sec2} 
Lucas' theorem and its variations}\hfill {\bf 5}
\vspace{2mm}

\indent \ref{subsec2.1} Lucas' theorem \dotfill{\bf 5}\\ 
\indent \ref{subsec2.2} Some consequences and extensions of Lucas' theorem 
\dotfill{\bf 7} 
\bigskip

\noindent{\bf \ref{sec3} 
Lucas type congruences for prime powers}\hfill {\bf 10}
\vspace{2mm}

\indent \ref{subsec3.1} Wolstenholme type congruences \dotfill{\bf 10} \\
\indent \ref{subsec3.2} Variations of Lucas' theorem modulo prime powers
  \dotfill{\bf 11} \\
\indent \ref{subsec3.3} Characterizations of Wolstenholme primes 
\dotfill{\bf 17} 

\bigskip

\noindent{\bf \ref{sec4} The Lucas property and the  $p$-Lucas property}
\hfill {\bf 18}
\vspace{2mm}

\indent \ref{subsec4.1} The Lucas property and the double Lucas property
 \dotfill{\bf 18} \\
\indent \ref{subsec4.2} Further Lucas type congruences  \dotfill{\bf 23} 

\bigskip

\noindent{\bf \ref{sec5} Lucas type theorems for some generalized binomial 
coefficients}\hfill {\bf 27}
\vspace{2mm}

\indent \ref{subsec5.1} Generalized binomial coefficients and related Lucas type congruences
 \dotfill{\bf 27} \\
\indent \ref{subsec5.2} Lucas type congruences for some classes of 
Lucas $u$-nomial coefficients   \dotfill{\bf 32} 

\bigskip

\noindent{\bf \ref{sec6} Some applications of Lucas' theorem}\hfill 
{\bf 36}
\vspace{2mm}

\indent \ref{subsec6.1} Lucas' theorem and the Pascal's triangle  \dotfill{\bf 36} \\
\indent \ref{subsec6.2} Another applications of Lucas's theorem  
\dotfill{\bf 40} 

\bigskip

\noindent{\bf References }\hfill {\bf 43}
\bigskip

 \noindent{\bf Appendix}\hfill {\bf 49}

\vfill\eject

\section{Introduction}\label{sec1}

Prime numbers have been studied since the earliest days of mathematics.
Congruences modulo primes have been widely investigated since the time
of Fermat. There are numerous useful and often remarkable congruences 
and divisibility results for binomial coefficients;
see \cite[Ch. XI]{d} for older results and \cite{gr} for a 
modern perspective.

 Let $p$ be a prime. Then by  {\it Fermat little theorem}, for each 
integer $a$ not divisible by $p$
  $$
 a^{p-1}\equiv 1\pmod{p}.
  $$
Furthermore, by {\it Wilson theorem, for any prime $p$}
  $$
(p-1)!+1\equiv 0\pmod {p}.
  $$

In attempting to discover some analogous expression 
which should be divisible by $n^2$, 
whenever $n$ is a prime, but not divisible if $n$ is a
composite number,  in 1819 Charles Babbage \cite{bab} is led 
to the congruence 
  \begin{equation*}
{2p -1\choose p-1}\equiv 1\pmod{p^2}
  \end{equation*}
for all primes $p\ge 3$. In 1862 J. Wolstenholme \cite{w}
proved that the above congruence holds modulo $p^3$ 
for any prime $p\ge 5$. 

 The study of arithmetic properties of binomial 
coefficients has a rich history.
As noticed in  \cite{gr}, many great mathematicians of the nineteenth century
considered problems involving binomial coefficients modulo a prime 
power (for instance Babbage \cite{bab}, Cauchy, Cayley, 
Gauss \cite{ga}, Hensel, Hermite \cite{he}, Kummer \cite{k}, 
Legendre, Lucas \cite{l} and \cite{l2}, and Stickelberger). 
They discovered a variety of elegant and surprising theorems which are
often easy to prove. For more information on these classical  results, their
extensions,  and new results about this subject, see books of 
Dickson \cite[Chapter IX]{d} and Guy \cite{gu}, while a more modern treatment 
of the subject is given by A.  Granville \cite{gr}.

Suppose that a prime $p$ and pair of integers  $n\ge m\ge 0$ 
are given. A beautiful {\it theorem of E. Kummer} of 1852 
(\cite[pp. 115--116]{k}; also see  \cite[p. 270]{d}) states that 
{\it the exact power of the prime $p$ which divides ${n\choose m}$ is 
given by the number of ``carries'' when $m$ and $n-m$ are added in 
base $p$ arithmetic}.
This is a fundamental result in the study of divisibility properties
of binomial coefficients.

 If $n=n_0+n_1p+\cdots +n_sp^s$ and
$m=m_0+m_1p+\cdots +m_sp^s$ are the $p$-{\it adic expansions} of 
nonnegative integers $n$ and $m$ 
(so that $0\le m_i,n_i\le p-1$ for each $i$),
then by {\it Lucas's theorem} 
established by \'{E}douard  Lucas in 1878 \cite{l} 
(also see \cite[p. 271]{d} and \cite{gr}), 
   $$
{n\choose m}\equiv \prod_{i=0}^{s}{n_i\choose m_i}\pmod{p}.
  $$
The same result is without proof also presented by Lucas  in 
1878, in Section XXI of  his massive journal paper \cite[pp. 229--230]{l2}.

This remarkable result by Lucas provides a simple way to compute
the binomial coefficient ${n\choose m}$ modulo a prime $p$ 
in terms of the binomial coefficients of the base-$p$ digits of 
$n$ and $m$. The above  congruence, the so-called  Lucas' theorem 
(or Theorem of Lucas) 
is a very important congruence in Combinatorial Number Theory and 
Combinatorics. In particular, this concerns the divisibility 
of binomial coefficients by primes. 
In this article, consisting of six sections, 
we provide a historical survey of Lucas type congruences, generalizations
of Lucas' theorem modulo prime powers and  Lucas like
theorems for some classes of generalized binomial coefficients.
Furthermore, we present some known applications 
of Lucas' theorem and certain of its variations in Number Theory and 
Combinatorics.

This  article is organized as follows. 
In Section \ref{sec2} we  mention several known algebraic and combinatorial 
proofs of Lucas' theorem. We also give  some  consequences 
and variations of Lucas' theorem. In Section \ref{sec3} we present 
a number of extensions and variations of Lucas' theorem modulo prime powers.  
In Section \ref{sec4} we consider the notions of the Lucas property and the double  
Lucas property. In this section we also present numerous integer sequences 
satisfying one of these properties or a certain similar Lucas type congruence. 
In particular, these properties are closely related to 
the divisibility properties of certain  binomial coefficients, matrices, 
different binomial sums, Ap\'{e}ry numbers, Delannoy numbers, Stirling 
numbers of the first and second kind etc. In Section \ref{sec5} we collect several known 
Lucas type congruences for some generalized binomial coefficients. In particular, this 
concerns the Fibonomial coefficients, the Lucas $u$-nomial coefficients, 
the Gaussian $q$-nomial coefficients and some their generalizations.
Finally, applications of Lucas' theorem are given in Section \ref{sec6}
of this survey article. 
Some of these applications are closely related to the determination
of number of entries of Pascal's triangle with a prescribed 
divisibility property. We also present some known primality criteria 
whose proofs are based on Lucas' theorem. Furthermore, 
we give certain known results concerning the characterizations 
of the algebraicity of some  classes of  formal power series
in terms of the notion of the $p$-Lucas property.      
      
\vfill\eject
\section{Lucas' theorem and its variations}\label{sec2}

\subsection{Lucas' theorem}\label{subsec2.1}

As noticed above, if $n=n_0+n_1p+\cdots +n_sp^s$ and
$m=m_0+m_1p+\cdots +m_sp^s$ are the $p$-{\it adic expansions} of integers 
$n$ and $m$  such that $0\le m_i,n_i\le p-1$ for each $i=0,1,\ldots,s$,
then a beautiful {\it Lucas's theorem} (\cite{l}; also see \cite{gr}
(\cite{l} and \cite[p. 271]{d}) states that 
 \begin{equation}\label{con1}
{n\choose m}\equiv \prod_{i=0}^{s}{n_i\choose m_i}\pmod{p}.
  \end{equation}
(with the usual convention that ${0\choose 0}=1$, and
${l \choose r}=0$ if $l<r$).
The  congruence \eqref{con1} was established by Lucas by  considering 
patterns in Pascal's triangle.
Furthermore,  \eqref{con1} is equivalent 
to the following  Lucas' earlier generalization \cite[p. 52]{l} of an 1869
result of H. Anton \cite[pp. 303--306]{an}
(also see \cite[p.  271]{d}):
 \begin{equation}\label{con2}
{n\choose m}\equiv {n\,{\rm div}\, p\choose m\,{\rm div}\,p}
{n\bmod{\,p}\choose m\bmod{p}}\pmod{p},
  \end{equation}
{\it where $n\,{\rm div}\, p$ denotes the integer quotient of $n$ by 
a prime $p$,
and $n\bmod{p}$ its remainder}. 
The congruence \eqref{con2} is in fact the equivalent form of Lucas' theorem
which is often stated in the follwing way:  
  \begin{equation}\label{con4}
{np +r\choose mp+s}\equiv {n \choose m}{r \choose s}\pmod{p},
  \end{equation}
where $p$ is a prime, $n,m,r$ and $s$ are nonnegative 
integers such that $0\le r,s\le p-1$. 
 
If a prime $p$ divides ${n\choose m}$ then \eqref{con1} follows easily from 
Kummer's theorem. However, if $p^l$ is the exact power of $p$ dividing ${n\choose m}$,
then we might ask for the value of $\frac{1}{p^l}{n\choose m}(\bmod{\,p})$.
The related result was discovered by H. Anton in 1869 
\cite{an} (see also \cite{gr}, \cite[pp.  3--4]{kaz} and \cite{sin})
who proved  {\it that if $p^l$ is the exact power of $p$ dividing ${n\choose m}$,
$($$l$ is by Kummer's theorem, the number of ``carries'' when $m$ and $n-m$ are 
added in base $p$ arithmetic$)$,  then
  \begin{equation}\label{con3}
\frac{(-1)^l}{p^l}{n\choose m}\equiv \frac{n_0!}{m_0!r_0!}\cdot
\frac{n_1!}{m_1!r_1!}\cdots\frac{n_s!}{m_s!r_s!}
\pmod{p},
  \end{equation}
where $n=n_0+n_1p+\cdots +n_sp^s$, 
$m=m_0+m_1p+\cdots +m_sp^s$, and   $r=n-m=r_0+r_1p+\cdots +r_sp^s$ 
with $0\le m_i,n_i,r_i\le p-1$ for each $i=0,1,\ldots,s$}.

\vspace{1mm}

{\it Remark} 1. Numerous authors have asked whether there is an analogous 
congruence modulo $p^l$ to \eqref{con3}, for arbitrary $l\ge 1$. 
In 1995 A. Granville \cite[Theorem 1]{gr} gave a positive answer to this 
question (see the congruence \eqref{con28}) in Subsection \ref{subsec3.2}). 
   \hfill $\Box$

\vspace{1mm}

 The several proofs offered for Lucas' theorem are 
primarily of to types-algebraic and combinatorial. 
The well known algebraic proof of Lucas' theorem due to N.J. Fine \cite{fi}
in 1947 is based on the binomial theorem for expansion of $(1+x)^n$. 
This proof runs as follows. Since by Kummer's theorem, 
the binomial coefficient ${p\choose k}$ is divisible by a prime $p$ for every 
$k=1,2,\ldots,p-1$, by the binomial expansion  it follows that 
 $$
(1+X)^p\equiv 1+X^p\pmod{p}.
  $$ 
Continuing by induction, we have that for every nonnegative integer $i$ 
  $$
(1+X)^{p^i}\equiv 1+X^{p^i}\pmod{p}.
  $$ 
Write $n$ and $m$ in base $p$, so that $n=\sum_{i=1}^sn_i$
and $m=\sum_{i=1}^sm_i$ for some nonnegative integers 
$s,n_0,\ldots,n_s,m_0,\ldots,m_s$  with $0\le n_i,m_i\le p-1$ 
for all $i=0,1,\ldots,s$. Then 
   \begin{equation*}\begin{split}
\sum_{m=0}^n{n\choose m}X^m &=(1+X)^n=
\prod_{i=0}^s\left((1+X)^{p^i}\right)^{n_i}\\
&\equiv \prod_{i=0}^s\left(1+X^{p^i}\right)^{n_i}=
\prod_{i=0}^s\left(\sum_{m_i=0}^{n_i}{n_i\choose m_i}X^{m_ip^i}\right)\pmod{p}\\
&= \prod_{i=0}^s\left(\sum_{m_i=0}^{p-1}{n_i\choose m_i}X^{m_ip^i}\right)\\
&=\sum_{m=0}^n\left( \prod_{i=0}^s{n_i\choose m_i}\right)X^m\pmod{p}.
 \end{split}\end{equation*}
By comparing the coefficients of $X^m$ on the left hand side 
and on the right hand side of the above congruence immediately
yields  Lucas' theorem given by \eqref{con1}.

As an application of a counting technique due to  M. Hausner in 1983 \cite{h},
in the same paper \cite[Example 4]{h} the author 
established another combinatorial proof of 
\eqref{con4}. Another proof of the congruence \eqref{con4} 
based on a simple combinatorial lemma is presented in 2005 by 
P.G. Anderson, A.T. Benjamin and J.A. Rouse 
in \cite[p. 268]{abr} (see also \cite{bq}). Another two proofs 
of Lucas' theorem, based on techniques from Elementary Number Theory
were obtained in 2010 by  S.-C. Liu and J.C.-C. Yeh \cite{ly}
and in 2012 by A. Laugier and M.P. Saikia \cite{ls2}.

The congruence \eqref{con4}  immediately yields 
 \begin{equation}\label{con5}
{np \choose mp}\equiv {n \choose m}\pmod{p}
  \end{equation}
since the same products of binomial coefficients are formed on 
the right side of Lucas's theorem in both cases, other than an  extra
${0\choose 0}=1$. 

\vspace{1mm}

A direct proof of the  congruence \eqref{con5}, 
 based on a polynomial method, is given in  \cite[Solution of Problem A-5, 
p. 173]{pu} as follows. It is well known that 
${p\choose i}\equiv 0(\bmod{\,p})$
for each $i=1,2,\ldots,p-1$ (see (\ref{con11})) or equivalently that
in the ring $\Bbb Z_{p}[x]$ we have $(1+x)^p=1+x^p$, where $\Bbb Z_p$ 
is the field of the integers modulo $p$. Thus in $\Bbb Z_{p}[x]$,
  $$
\sum_{k=0}^{np}{np\choose k}x^k=(1+x)^{np}=\left((1+x)^p\right)^n=
(1+x^p)^n=\sum_{j=0}^n{n\choose j}x^{jp}.
  $$
Since coefficients of like powers must be congruent modulo $p$ in the equality
  $$
\sum_{k=0}^{np}{np\choose k}x^k=\sum_{j=0}^n{n\choose j}x^{jp}
  $$
in $\Bbb Z_{p}[x]$, we see that 
   $$
{np\choose mp}\equiv {n\choose m}\pmod{p}\quad for\,\, m=0,1,\ldots,n.
  $$

Further, notice that the Lucas' congruence \eqref{con4} 
easily follows by induction on the sum $r+s\ge 0$ using 
the base induction $r+s=0$ with $r=s=0$ satisfying via the  
congruence \eqref{con5}, and the Pascal formulas:
  $$
{np+(r+1)\choose mp+s}={np+r\choose mp+(s-1)}+{np+r\choose mp+s}
 $$
and
 $$
{np+r \choose mp+(s+1)}={np+(r-1)\choose mp+s}+{np+(r-1)\choose mp+(s+1)}.
 $$

\vspace{1mm}

{\it Remark} 2. The Lucas' congruence (\ref{con4}) also can be interpreted 
as a result about cellular automata (cf. Granville \cite[Section 5]{gr}). 
Namely, Lucas' theorem can be interpreted as a two-dimensional 
$p$-automaton (for a formal definition see  \cite{ahps}). 
\hfill$\Box$

\subsection{Some consequences and extensions of Lucas' theorem}\label{subsec2.2}

Here, as always in the sequel, $p$ will denote any prime.

As noticed in 2011 by A. Nowicki \cite[the congruences 7.3.1--7.33]{no}, 
{\it if $n=n_0+n_1p+\cdots +n_sp^s$ is the $p$-adic expansion 
of a positive integer $n$, then for each $k=0,1,\ldots,s$
 \begin{equation}\label{con6}
{n \choose p^k}\equiv n_k\equiv \left\lfloor\frac{n}{p^k}\right\rfloor\pmod{p},
  \end{equation}
holds}, and consequently, 
  \begin{equation}\label{con7}
{n \choose p}\equiv \left\lfloor\frac{n}{p}\right\rfloor\pmod{p},
  \end{equation}
{\it where $\lfloor x \rfloor$ is  the greatest   integer less than or equal 
to $x$}. 
\vspace{1mm}

{\it Remark} 3. The congruence \eqref{con7} 
is  proposed by L.E. Clarke \cite{cl} in 1956 as a problem which is 
solved in 1957 by P.A. Piza \cite{pi}.\hfill $\Box$ 
\vspace{1mm}

Moreover, {\it if $0\le r<p^f$ and $0\le m<p^f$, then 
the Lucas' congruence \eqref{con4} immediately yields $($see 
\cite[the congruence 7.3.6]{no}$)$}
 \begin{equation}\label{con8}
{p^f+r \choose m}\equiv {r\choose m}\pmod{p}.
  \end{equation}    
Furthermore, {\it if $0\le r<p^f$, $0\le m<p^f$ and $a\ge 0$, then by 
Lucas' theorem $($see \cite[the congruence 7.3.7]{no}}),   
 \begin{equation}\label{con9}
{ap^f+r \choose m}\equiv {r\choose m}\pmod{p}.
  \end{equation}
Moreover, {\it if $0\le r<p^f$ and $p^f\le m$, then by 
\cite[the congruence 7.3.8]{no}},    
 \begin{equation}\label{con10}
{p^f+r \choose m}\equiv {r\choose m-p^f}\pmod{p}.
  \end{equation}

Lucas' theorem immediately yields the following 
 well known congruence:  
  \begin{equation}\label{con11}
{p\choose k} \equiv 0 \pmod{p},
   \end{equation}
where  $p$ is a prime and $k$ is an integer such that  $1\le k\le p-1$.

Furthermore, {\it if $p$ is a prime and $f$ a positive integer, then  
by Lucas' theorem  for any $f\ge 1$ and $1\le k\le p^f-1$ we 
have $($see, e.g., \cite[Theorem 24]{bq}$)$}
  \begin{equation}\label{con12}
{p^f\choose k} \equiv 0 \pmod{p}.
   \end{equation}
Further, {\it if $p$ is a prime and $n$, $m$ and $k$ are positive integers 
with $m\le n$, then the congruence \eqref{con5} by induction easily 
yields $($see \cite[Lemma 2.1]{me3}$)$ 
    \begin{equation}\label{con13} 
{np^k\choose mp^k} \equiv {n\choose m} \pmod{p}.
    \end{equation}} 

An alternative version of Lucas' theorem 
was noticed in 1994 by J. M. Holte \cite[p. 60]{hol2} 
(also  see  \cite[p. 227]{hol1})  as follows: {\it If  
 $$
B(m,n):={m+n\choose m}=\frac{(m+n)!}{m!n!},
 $$
then 
  \begin{equation}\label{con52}
B(m,n)\equiv B(m\,{\rm div}\, p,n\,{\rm div}\, p)B(m\bmod{\,p},n\bmod{\,p}) 
\pmod{p},
   \end{equation}
where $m\,{\rm div}\, p$ is the integer quotient of $m$ by $p$
and $m\bmod{\,p}$ is the remainder of $m$ by division by $p$.
$($similarly, for $n$ instead of $m$$)$. 
It follows that {\it if $n=n_0+n_1p+\cdots +n_sp^s$
and $m=m_0+m_1p+\cdots +m_sp^s$, where $0\le m_i,n_i\le p-1$
for each $i=0,1,\ldots s$, then
 \begin{equation}\label{con53}
B(m,n)\equiv \prod_{i=0}^sB(m_i, n_i) \pmod{p}.
   \end{equation}}

Consequently,  {\it $p\mid B(m,n)$ if and only if 
 $p\mid B(m_i,n_i)$ for some $i\in\{0,1,\ldots,s\}$}}.

Following Granville \cite[Section 6]{gr}, for an integer polynomial $f(X)$ 
of degree $d$, define the numbers ${m\choose n}_f$ with $m,n\in \Bbb Z$ 
by the generating function 
   $$
f(X)^m=\sum_{n=0}^{md}{m\choose n}_fX^n,
   $$
and let ${m\choose n}_f=0$ if $n<0$ or $n>md$ (note that 
${m\choose n}_f={m\choose n}$ when $f(X)=X+1$). Clearly, 
by Fermat little theorem, $f(X)^p\equiv f(X^p)(\bmod{\, p})$, and 
using this in 1995 A. Granville \cite[Section 6, the congruence (24)]{gr} 
proved the following generalization of the congruence \eqref{con3}: 
{\it If $p$ is a prime, $m,n$  nonnegative integers such that $m=pl+m_0$,
 $n=pt+n_0$,  $l,t,m_0,n_0\in\Bbb N$ and $0\le m_0,n_0\le p-1$,
then
     \begin{equation}\label{con14} 
{m\choose n}_f \equiv \sum_{k=0}^{d-1}{\lfloor m/p\rfloor \choose 
\lfloor n/p\rfloor -k}_f{m_0\choose n_0+kp}_f \pmod{p}.
    \end{equation}}
Notice that when $f(X)=X+1$ then the congruence \eqref{con14}  becomes 
    $$
{m\choose n}\equiv {\lfloor m/p\rfloor \choose 
\lfloor n/p\rfloor}{m_0\choose n_0}\pmod{p},
   $$
which is in fact the Lucas's congruence \eqref{con4}.

By using a congruence based on Burnside's theorem, 
in 2005, T.J. Evans \cite[Theorem 3]{ev} 
proved the following extension of Lucas' theorem 
involving Euler's totient function $\varphi$: 
{\it If $n\ge 1$, $m,M,m_0,r,R$ and $r_0$ are nonnegative integers 
such that  $m=Mn+m_0$, $r=Rn+r_0$, with $0\le m_0, r_0<n$, then 
  \begin{equation}\label{con15+}
\sum_{d\mid n}\varphi\left(\frac{n}{d}\right)
\sum_{j=-(d-1)}^{d-1}\sum_{\lVert a\rVert_d=\atop R-(j/d)}{M\choose a_1}\cdots
{M\choose a_d}{m_0\choose r_0+(n/d)j} \equiv 0\pmod{n}, 
 \end{equation}
where the summation runs among all positive divisors $d$ of $n$.}

\vspace{1mm}
{\it Remark} 4.
It was proved in \cite[Corollary 3]{ev} 
that Lucas' theorem easily follows from the congruence \eqref{con15+}. \hfill $\Box$

\section{Lucas type congruences for prime powers}\label{sec3}

\subsection{Wolstenholme type congruences}\label{subsec3.1}
Notice that for any prime $p$ the  congruence \eqref{con5} with $n=2$ and 
$m=1$ becomes 
   $$
{2p \choose p}\equiv 2\pmod{p},
  $$
whence by the identity ${2p\choose p}=2{2p-1\choose p-1}$
it follows that {\it for any prime $p$}
  \begin{equation}\label{con16}
{2p -1\choose p-1}\equiv 1\pmod{p}.
  \end{equation}
As noticed in \ref{sec1}, in 1819 Charles Babbage  \cite{bab} (also 
see \cite[Introduction]{gr} or \cite[page 271]{d})
showed that the congruence (\ref{con16}) holds modulo $p^2$,
that is, {\it for a prime $p\ge 3$ holds}
  \begin{equation}\label{con17}
{2p -1\choose p-1}\equiv 1\pmod{p^2}.
  \end{equation}
\vspace{1mm}

{\it Remark} 5.  A combinatorial proof 
of the congruence \eqref{con17} can be found in \cite[Exercise 14(c) on 
page 118]{st}.
\hfill $\Box$ 
\vspace{1mm}

The congruence \eqref{con17} was generalized in 1862 by Joseph Wolstenholme  
\cite{w} as it is presented in the next section.  
Namely, {\it Wolstenholme's theorem} asserts that   
  \begin{equation}\label{con18}
{2p -1\choose p-1}\equiv 1\pmod{p^3}
  \end{equation}
for all primes $p\ge 5$. 

For a survey of Wolstenholme's theorem see \cite{me5} and for its 
extensions  see \cite{z1} and \cite{me10}.

By Glaisher's congruence  
\cite[p. 323]{gl2} (also see  \cite[Section 6]{me5}),   
{\it for any positive integer $n$ and a prime $p\ge 5$ holds} 
    \begin{equation*}
{np-1\choose p-1}\equiv 1\pmod{p^3},
      \end{equation*}
which by the identity ${np\choose p}=n{np-1\choose p-1}$
yields \cite[the congruence 7.1.5]{no}  
  \begin{equation}\label{con19}
{np\choose p}\equiv n\pmod{p^3}.
      \end{equation}

In 1949 W. Ljunggren \cite{bs}  generalized the congruence \eqref{con19} 
as follows (also see \cite[Theorem 4]{ba}, \cite{gr} and 
\cite[Problem 1.6 (d)]{st}, and for a simple proof see \cite{sio}):
{\it if $p\ge 5$ is a prime, $n$ and $m$ are positive integers
with $m\le n$, then}
   \begin{equation}\label{con20}
{np\choose mp} \equiv {n\choose m} \pmod{p^3}.
   \end{equation}

\vspace{1mm}

{\it Remark} 6. Ljunggren's congruence \eqref{con20} is refined modulo 
$p^5$ in 2007 by J. Zhao \cite[Theorem 3.5]{z2}.\hfill$\Box$
\vspace{1mm}

{\it Remark} 7.  Note that the  congruence \eqref{con20} with $m=1$ and 
$n=2$ reduces to the Wolstenholme's congruence 
\eqref{con18}.  \hfill $\Box$.
\vspace{1mm}

Further, the congruence \eqref{con20}
is refined in 1952 by E. Jacobsthal  \cite{bs} (also see \cite{gr})
 as follows: {\it if $p\ge 5$ is a prime, $n$ and $m$ are positive integers
with $m\le n$, then
    \begin{equation}\label{con21}
{np\choose mp} \equiv {n\choose m} \pmod{p^t},
   \end{equation}
where $t$ is the power of $p$ dividing 
$p^3nm(n-m)$ $($this exponent $t$ can only be 
increased if $p$ divides  $B_{p-3}$, the $(p-3)$rd  Bernoulli number$)$}.

\vspace{1mm}
{\it Remark} 8. In the literature, the congruence 
\eqref{con21} is often called {\it Jacobsthal-Kazandzidis congruence}
(see e.g., \cite[Section 11.6, p. 380]{co}).  \hfill $\Box$
  \vspace{1mm} 

In 2008 C. Helou and G. Terjanian \cite[the congruence (1) of Corollary on 
page 490]{ht} refined the Jacobsthal's result as follows 
(also see \cite[Section 11.6, Corollary 11.6.22, p. 381]{co} for 
a stronger form)): {\it If $p\ge 5$ is a prime, $n$ and $m$ are positive 
integers with $m\le n$, then
    \begin{equation}\label{con22}
{np\choose mp} \equiv {n\choose m} \pmod{p^t},
   \end{equation}
where $t$ is the power of $p$ dividing $p^3m(n-m){n\choose m}$.}

By a problem N4 of Short list of 48th IMO 2006 \cite{djp}, 
{\it for every integer $k\ge 2$, $2^{3k}$ divides the number
  \begin{equation}\label{con22+}  
{2^{k+1}\choose 2^k}-{2^k\choose 2^{k-1}}
  \end{equation}
but $2^{3k+1}$ does not}.

\subsection{Variations of Lucas' theorem modulo prime powers}\label{subsec3.2}

In 1991 D.F. Bailey \cite[Theorem 4]{ba2} proved that
{\it if $p$ is a prime, $n$ and $r$ are  nonnegative integers 
and $s$ a positive integer less than $p$, then}
 \begin{equation}\label{con23}
{np\choose rp+s}\equiv (r+1){n\choose r+1}{p\choose s}\pmod{p^2}.
      \end{equation}
In the same paper \cite[Theorem 5]{ba2}, the author
extended the previous congruence as follows:   
{\it if $p\ge 5$ is a prime, $0\le m\le n$, $0\le r<p$ 
and $1\le s<p$, then}
 \begin{equation}\label{con24}
{np^2\choose mp^2+rp+s}\equiv (m+1){n\choose m+1}{p^2\choose rp+s}\pmod{p^3}.
      \end{equation}

\vspace{1mm}
{\it Remark} 9. Notice  that Bailey's  proof of the congruence \eqref{con24}
(proof of Theorem 5 in \cite{ba}) is deduced  applying  the Ljunggren's 
congruence \eqref{con20} (Theorem 4 in \cite{ba}) and a counting technique of 
M. Hausner from  \cite{h}. \hfill$\Box$
 \vspace{1mm}

In 1992 D.F. Bailey \cite[Theorem 2.1]{ba3}
generalized his congruence \eqref{con24} modulo any prime power as follows:   
{\it if $p\ge 5$ is a prime, $0\le m\le n$, 
$s\ge 1$, and $a_0,a_1, \ldots ,a_{s-1}$ are nonnegative integers
such that $1\le a_0<p$ and $0\le a_k<p$ for every $k=1,2,\ldots, s-1$,
then}
 \begin{equation}\label{con25}\begin{split}
&{np^s\choose mp^s+a_{s-1}p^{s-1}+\cdots +a_1p+a_0}\\
&\equiv (m+1){n\choose m+1}{p^s\choose a_{s-1}p^{s-1}+\cdots +a_1p+a_0}\pmod{p^{s+1}}.
      \end{split}\end{equation}

{\it Remark} 10. If we put $a=a_{s-1}p^{s-1}+\cdots +a_1p+a_0$,
then the congruence \eqref{con25} can be written as 
  \begin{equation}\label{con26}
{np^s\choose mp^s+a}\equiv (m+1){n\choose m+1}{p^s\choose a}\pmod{p^{s+1}},
      \end{equation}
where $a$ is a positive integer less than $p^s$ which is not divisible by $p$. 
\hfill $\Box$.
\vspace{1mm}

Using a multiple application of Lucas' theorem, 
in  2012  the author  of this article  \cite[Theorem 1.1]{me8} 
proved the following similar congruence to \eqref{con26}:
    \begin{equation}\label{con27}
{np^s\choose mp^s+a}\equiv (-1)^{a-1}a^{-1}(m+1){n\choose m+1}p^s
\pmod{p^{s+1}},
      \end{equation}
{\it where $p$ is a prime, $n$, $m$, $s$ and $a$ are  nonnegative 
integers such that $n\ge m$, $s\ge 1$, $1\le a\le p^s-1$, 
and $a$ is not divisible by $p$}.

\vspace{1mm}

{\it Remark} 11. The congruence \eqref{con26} is an immediate  consequence
of the congruence \eqref{con27} (see \cite[Corollary 1.2 and its 
proof]{me8}).\hfill$\Box$

\vspace{1mm}

In 1990 D.F. Bailey \cite[Theorem 3]{ba} 
(cf.  \cite[Theorem with $k=2$]{me4}) proved  the following result:   
{\it If $p$ is a prime, $n,m,n_0$ and $m_0$  are
 nonnegative integers, and  $n_0$ and $m_0$ are both less than $p$, then}
  \begin{equation}\label{con38}
{np^2 +n_0 \choose mp^2+m_0}\equiv {n\choose m}
{n_0 \choose m_0}
\pmod{p^2}.
  \end{equation}

Furthermore, in the same paper Bailey \cite[Theorem 5]{ba} 
(cf. \cite[Theorem with $k=3$]{me4}) extended the above result as 
 follows: {\it If $p$ is a prime greater than $3$
and $n,m,n_0$ and $m_0$  are nonnegative integers such that $n_0$ and $m_0$ 
are less than $p$, then}
  \begin{equation}\label{con38+}
{np^3 +n_0 \choose mp^3+m_0}\equiv {n\choose m}
{n_0 \choose m_0}
\pmod{p^3}.
  \end{equation}

Kummer's theorem given in Section \ref{sec1}, is useful 
in situations where the binomial coefficient is divisible
by a prime power. However, if the binomial coefficient is not
congruent to zero modulo a prime,  then the question remains for a way to 
simplify the expression.     
In 1995 A. Granville \cite[Theorem 1]{gr} generalized Anton's congruence 
\eqref{con3} modulo prime powers as follows.  
{\it For a given integer $k$ define 
$(k!)_p$ to be the product of all integers  less than or equal to 
$k$, which are not divisible by $p$. Suppose that prime power $p^f$ and 
positive integers $n$ and $m$ are given with $r:=n-m\ge 0$. Write  
$n=n_0+n_1p+\cdots +n_sp^s$ in base $p$, and let $N_j$ be the least positive 
residue of $\lfloor n/p^j\rfloor (\bmod{\, p^f})$ for each $j\ge 0$   
$($so that  $N_j=n_j+n_{j+1}p +\cdots +n_{j+f-1}p^{f-1}$$);$ also make 
the corresponding definitions for $m_j,M_j,r_j,R_j$. Let $e_j$ be the number 
of indices $i\ge j$ for which $n_i<m_i$ $($that is, the number of ``carries''
when adding $m$ and $r$ in base $p$, on or beyond the $j$th digit$)$. Then 
 \begin{equation}\label{con28}
\frac{1}{p^{e_0}}{n\choose m}\equiv (\pm 1)^{e_{f-1}}
\frac{(N_0!)_p}{(M_0!)_p(R_0!)_p}\cdot \frac{(N_1!)_p}{(M_1!)_p(R_1!)_p}
\cdots \frac{(N_s!)_p}{(M_s!)_p(R_s!)_p}\pmod{p^f},
  \end{equation}
where  $(\pm 1)$ is $(-1)$ except if $p=2$ and $f\ge 3$}.

Here, as usually in the sequel, we will consider the  congruence relation 
modulo a prime power $p^l$ extended to the ring of rational numbers
with denominators not divisible by $p$. 
For such fractions we put $m/n\equiv r/s \,(\bmod{\,\,p^l})$ 
if and only if $ms\equiv nr\,(\bmod{\,\,p^l})$, and the residue
class of $m/n$ is the residue class of $mn'$ where 
$n'$ is the inverse of $n$ modulo $p^l$.  

A result which gives readily an extension of 
Lucas' theorem in the form of the  congruence 
to prime power moduli is given in 1992 by A. Granville 
\cite[Proposition 2]{gr1} as follows: {\it 
For each positive integer $j$, define  $n_j$ to be  the least nonnegative 
residue of an integer $n$ modulo $p^j$. If $p$ is a prime that does not 
divide ${n\choose m}$, then  
   \begin{equation}\label{con29} 
{n\choose m} 
\equiv {\lfloor n/p\rfloor\choose \lfloor m/p\rfloor}{n_f\choose m_f}
\Bigg/{\lfloor n_f/p\rfloor\choose \lfloor m_f/p\rfloor} \pmod{p^f},
    \end{equation} 
for any positive integer $f$}.

 In particular, {\it if  ${n\choose m}$  is not 
divisible by $p$ and $m\equiv n (\bmod{\,p^f})$, then  
by \eqref{con29} $($also see \cite[the congruence 7.1.16]{no}$)$}   
   \begin{equation}\label{con30}
{n \choose m}\equiv {\left\lfloor n/p\right\rfloor
\choose \left\lfloor m/p\right\rfloor}\pmod{p^f}.
  \end{equation}

As observed in 1998 by D. Berend and J.E. Harmse 
\cite[p. 34, congruence (2.2)]{bh}, 
{\it if a prime $p$ does not divide ${n\choose m}$ and
$n=n_0+n_1p+\cdots +n_sp^s$, 
$m=m_0+m_1p+\cdots +m_sp^s$ are the $p$-adic expansions of $n$ and $m$, then
 iterating the  congruence \eqref{con29}, we find that
    \begin{equation}\label{con31} 
{n\choose m}\equiv \frac{P}{Q}\pmod{p^f}, 
 \end{equation} 
where 
  $$
P=\prod_{i=0}^{k-f+1}{n_i+n_{i+1}p+\ldots +n_{i+f-1}p^{f-1}\choose
m_i+m_{i+1}p+\ldots +m_{i+f-1}p^{f-1}}
  $$
and} 
 $$
Q=\prod_{i=1}^{k-f+1}{n_i+n_{i+1}p+\ldots +n_{i+f-2}p^{f-2}\choose
m_i+m_{i+1}p+\ldots +m_{i+f-2}p^{f-2}}.
  $$
The congruence \eqref{con31} was established in 1991 independently 
by K. Davis and W. Webb \cite[Theorem 3]{dw1} 
(also see  \cite[p. 88, Theorem 5.1.2]{lo}), which is there 
formulated as follows: {\it If $n=n_0+n_1p+\cdots +n_sp^s$, 
$m=m_0+m_1p+\cdots +m_sp^s$ are the $p$-adic expansions of $n$ and $m$, 
and $l<s$, then
   \begin{equation}\label{con32}\begin{split}
{n\choose m} \equiv &  
{n_0+n_1p+\cdots +n_{s-1}p^{s-1}+n_sp^s
\choose m_0+m_1p+\cdots +m_{s-1}p^{s-1}+m_sp^s}\\
\equiv & 
\frac{{n_{s-l}+\cdots +n_sp^{s-l}\choose m_{s-l}+\cdots +m_sp^{s-l}}
\cdots {n_0+\cdots +n_lp^{s-l}\choose m_{0}+\cdots +m_lp^{s-l}}}
{{n_{s-l+1}+\cdots +n_{s-1}p^{s-l-1}\choose m_{s-l+1}+\cdots +
m_{s-1}p^{s-l-1}}\cdots 
{n_0+\cdots +n_{l-1}p^{s-l-1}\choose m_0+\cdots +m_{l-1}p^{s-l-1}}}
   \pmod{p^l}.
 \end{split}\end{equation}
If $a=a_0+a_1p+\cdots +a_{k-1}p^{k-1}+a_kp^k$ and 
$b=b_0+b_1p+\cdots +b_{k-1}p^{k-1}+b_kp^k$ are the $p$-adic expansions 
of $a$ and $b$ such that $b_k>a_k$, then we define 
   $$
{a_0+a_1p+\cdots +a_{k-1}p^{k-1}+a_kp^k  \choose 
b_0+b_1p+\cdots +b_{k-1}p^{k-1}+b_kp^k} =
p{a_0+a_1p+\cdots +a_{k-1}p^{k-1}  \choose 
b_0+b_1p+\cdots +b_{k-1}p^{k-1}}.
  $$}
\vspace{1mm}

{\it Remark} 12.  For help in understanding the above result concerning 
the congruence \eqref{con32}, we offer the following example 
\cite[p. 88]{lo}:
   \begin{equation*}\begin{split}
{386\choose 154}&={3\cdot 11^2+2\cdot 11+1\choose 11^2+3\cdot 11}
\equiv\frac{{3\cdot 11+2\choose 11+3}{2\cdot 11+1\choose 3\cdot 11}}{{2\choose 3}}\pmod{11^2}\\
\qquad \qquad &\equiv {3\cdot 11+2\choose 11+3}{1\choose 0}\equiv {35\choose 14}\pmod{11^2}.  
\qquad\qquad \qquad\hfill \Box
   \end{split}\end{equation*}
\vspace{1mm}

In 2005  A.D. Loveless \cite[p. 88]{lo} noticed that the above result 
concerning  the congruence \eqref{con32} can be used to simplify 
general classes of congruences modulo prime powers  involving
binomial coefficients. In particular, 
Loveless \cite[p. 88, Theorem 5.1.3]{lo}) proved that {\it if 
$p$ is a prime, $s$ and $n$ are  positive integers with $n\le p^s$,  
then}
    \begin{equation}\label{con33}   
{p^s\choose n}\equiv \left\{
    \begin{array}{ll}
0 & (\bmod{\,p^s})\quad if\,\,n\not\equiv 0\,(\bmod{\,p})\\
{p^{s-1}\choose n/p}& (\bmod{\,p^s})\quad if \,\, p\mid n.
  \end{array}\right.
  \end{equation}

A similar result was earlier directly proved in 1980 by P.W. Haggard 
and J.O. Kiltinen \cite[p. 398, Theorem]{hk}. 
This result asserts that {\it if $p$ is a prime, $l$ and $f$ 
are positive integers with $f\ge l-1$ and $0\le n\le p^f$, then}
     \begin{equation}\label{con34}   
{p^f\choose n}\equiv \left\{
    \begin{array}{ll}
0 & (\bmod{\,p^l})\quad if\,\,n\not\equiv 0\,(\bmod{\,p^{f-l+1}})\\
{p^{l-1}\choose i}& (\bmod{\,p^l})\quad if \,\, n=i\cdot p^{f-l+1}.
  \end{array}\right.
  \end{equation}
Using the congruence \eqref{con32}, 
in 1993 K. Davis and W. Webb  \cite{dw2} generalized  Bailey's  results 
concerning the congruences  \eqref{con38} and  \eqref{con38+}
for any modulus $p^k$ with $p\ge 5$ and $k\ge 1$.
They proved \cite[Theorem 3]{dw2} that {\it if
 $p$ is any prime,  $k,n,m,a,b$  and $s$ 
are positive  integers such that $0<a,b<p^s$, then}
   \begin{equation}\label{con35}
{np^{k+s} +a \choose mp^{k+s}+b}\equiv {np^k\choose mp^k}
{a \choose b}\pmod{p^{k+1}}.
    \end{equation}
\vspace{1mm}

{\it Remark} 13. Notice that under the same assumption preceding the 
congruence \eqref{con35}, and if 
${np^{k+s} +a \choose mp^{k+s}}\not\equiv 0(\bmod{\, p})$,
then  the  congruence \eqref{con35} can be obtained  by 
iterating $s$ times the  Granville's congruence \eqref{con29}. 
Notice also that the condition 
${np^{k+s} +a \choose mp^{k+s}}\not\equiv 0(\bmod{\, p})$ is by Lucas' 
theorem equivalent to the following two conditions: 
${n \choose m}\not\equiv 0(\bmod{\, p})$ and 
${a \choose b}\not\equiv 0(\bmod{\, p})$.\hfill $\Box$ 
\vspace{1mm}

Further, by repeated application of the congruence \eqref{con35}, 
and using Ljunggren's congruence \eqref{con20}, we find that
{\it under the same assumptions preceding the congruence \eqref{con35}
\cite[Corollary 1]{dw2}  for any prime $p>3$,
   \begin{equation}\label{con36}
{np^{k+s} +a \choose mp^{k+s}+b}\equiv {np^{\lfloor k/3\rfloor}
\choose mp^{\lfloor k/3\rfloor}}{a \choose b}
\pmod{p^{k+1}}.
    \end{equation}}
 
In particular, the congruence \eqref{con36} with $s=1$ and $k-1\ge 0$ 
instead of $k$ implies that {\it for each prime $p\ge 5$
and for all integers $k\ge 1$, $n\ge 0$, $a$ and $b$ with  $0\le a, b<p$
   \begin{equation}\label{con37}
{np^k +a \choose mp^k+b}\equiv {np^{\lfloor(k-1)/3\rfloor}
\choose mp^{\lfloor(k-1)/3\rfloor}}{a \choose b}
\pmod{p^k}.
    \end{equation}}
Furthermore,  {\it  the congruence \eqref{con37} 
with $\lfloor k/2\rfloor$ instead of $\lfloor(k-1)/3\rfloor$
is satisfied for $p=2$, and   the congruence \eqref{con37}  
with $\lfloor (k-1)/2\rfloor$ instead of $\lfloor(k-1)/3\rfloor$
is also satisfied for $p=3$}.

\vspace{1mm}
{\it Remark} 14. As noticed above, a proof of the congruence 
\eqref{con36} given by Davis and Webb is based on their earlier result 
from \cite{dw1} given by the congruence \eqref{con36}. However, this result 
together with related proof is slightly more complicated. 
In 2012 the author of this article \cite[Theorem]{me4} 
gave a simple induction  proof of the congruence \eqref{con37} 
which  uses only the usual properties of 
binomial coefficients.  \hfill$\Box$
 \vspace{1mm}

Adapting  Fine's method \cite{fi}, in 1988 R.A. Macleod \cite[Theorem 2]{ma}
proved the following variation of Lucas' theorem:
{\it Let $p$ be a prime, let $r$ be a positive integer, and let
  $$
M=\sum_{i=0}^kM_ip^{ir}, \quad with \quad 0\le M_i<p^r \quad for\,\, all
\,\, i=0,1,\ldots,k.
  $$
Then for every nonnegative integer $N$ such that $0\le N\le M$
      \begin{equation}\label{con40}
{M\choose N}\equiv \sum {p^{r-1}M_0\choose N_0}{p^{r-1}M_1\choose N_1}
\cdots {p^{r-1}M_k\choose N_k}\pmod{p^r},
      \end{equation}
where the summation ranges over all $k+1$-tuples $(N_0,N_1,\ldots,N_k)$
such that 
 $$
p^{r-1}N=\sum_{i=0}^kN_ip^{ir}, \quad with \quad 0\le N_i<p^{r-1}M_i\quad
 for\,\, all \,\, i=0,1,\ldots,k.
  $$}

Quite recently, in 2014 E. Rowland and R. Yassawi 
\cite[Section 5, Theorem 5.3]{ro2} established a new generalization 
of Lucas' theorem to prime powers as follows: {\it Let $p$ be a prime, let $f$ be 
a positive integer and let $D=\{ 0,1,\ldots, p^f-p^{f-1}\}$.
 If $n=n_0+n_1p+\cdots +n_sp^s$ and 
$m=m_0+m_1p+\cdots +m_sp^s$ are the $p$-adic expansions of 
nonnegative integers $n$ and $m$, then
      \begin{equation}\label{con41}\begin{split}   
{n\choose m}\equiv &\sum_{(i_0,\ldots,i_l)\in D^{l+1}\atop 
(j_0,\ldots,j_l)\in D^{l+1}}(-1)^{n-i+\sum_{h=0}^li_h}{p^{f-1}-1\choose n-i}
{n-i\choose m-j}\\
&\times \prod_{h=0}^l{p^f-p^{f-1}\choose i_h}{i_h\choose j_h}\pmod{p^f},
     \end{split}\end{equation}
where $i=\sum_{h=0}^li_hp^h$ and $j=\sum_{h=0}^lj_hp^h$}.

\vspace{1mm}
{\it Remark} 15. Note that $i=\sum_{h=0}^li_hp^h$ and $j=\sum_{h=0}^lj_hp^h$ are 
representations of integers $i$ and $j$ in base $p$ 
with an enlarged digit set $D$ rather than the standard digit set
$\{0,1,\ldots,p-1\}$.\hfill$\Box$

\vspace{1mm}

{\it Remark} 16. E. Rowland and R. Yassawi \cite[Section 5]{ro2} showed that 
a broad range of multidimensional sequences possess ``Lucas products'' modulo 
a prime $p$. Furthermore, in 2009 K. Samol and D. van Straten 
\cite[Proposition 4.1]{ss} established the Lucas type congruence  
for a sequence  whose terms are constant terms of $P(x)^n$ for certain 
Laurent polynomials $P(x)$.\hfill $\Box$

\subsection{Characterizations of Wolstenholme primes}\label{subsec3.3}
A prime $p$ is said to be a {\it Wolstenholme prime} if it 
satisfies the congruence 
$$
{2p-1\choose p-1} \equiv 1 \, (\bmod{\, p^4}),
  $$
or equivalently,
 \begin{equation}\label{con43}
{2p\choose p} \equiv 2 \pmod{p^4}.
  \end{equation}
The two known such primes are 16843 and 2124679, and 
R.J. McIntosh and E.L. Roettger reported in \cite{mr} that these
primes are only two  Wolstenholme primes less than $10^9$.
However, McIntosh in \cite{m} conjectured
that there are infinitely many  Wolstenholme primes (for more information
see \cite{me6}).   
By the well known result of J.W.L. Glaisher in 1900 \cite[p. 323]{gl2} 
(also see \cite[the congruence (1.2)]{me1}), 
     \begin{equation}\label{con43}
{2p-1\choose p-1} \equiv 1-\frac{2}{3}p^3B_{p-3}\pmod{p^4},
    \end{equation}
where  $B_k$ ($k=0,1,2,\ldots$) are Bernoulli numbers  
defined by the generating function \cite{ir}
   $$
\sum_{k=0}^{\infty}B_k\frac{x^k}{k!}=\frac{x}{e^x-1}\,.
  $$
The  congruence \eqref{con43} shows that a prime $p$ is a Wolstenholme prime
if and only if  $p$ divides the numerator of $B_{p-3}$, 
the $(p-3)$rd  Bernoulli number.

As an application of the  congruences \eqref{con37} with $k=4$ and 
Jacobsthal's congruence \eqref{con21}, 
we can obtain the following characterization of  Wolstenholme primes
given in 2012 by the author of this article  \cite[Proposition]{me4}:
{\it The following statements about a prime
$p\ge 5$ are equivalent.
    \begin{itemize}
\item[(i)] $p$ is a Wolstenholme prime;
\item[(ii)] for  all  nonnegative integers $n$ and $m$ the congruence 
   \begin{equation}\label{con45}
{np\choose mp} \equiv {n\choose m} \pmod{p^4}
    \end{equation}
holds;
\item[(iii)]  for all nonnegative integers $n,m,n_0$ and $m_0$  such
that $n_0$ and $m_0$ are less than $p$,
   \begin{equation}\label{con46}
{np^4 +n_0 \choose mp^4+m_0}\equiv {n\choose m}
{n_0 \choose m_0}\pmod{p^4}.
  \end{equation}
   \end{itemize}}

\section{The Lucas property and the  $p$-Lucas property}\label{sec4}

\subsection{The Lucas property and the double Lucas property}\label{subsec4.1}

In 1992  R.J. McIntosh \cite{m2} proposed the following definition:

  \vspace{1mm}
{\it Definition}. The integer  sequence  $(a_n)_{n\ge 0}$ has the
{\it  Lucas property} if $a_0=1$, and for every prime $p$, every 
$n\ge 0$, and every $j\in \{0,1,\ldots,p-1\}$ the congruence
     \begin{equation}\label{con47}
a_{pn+j}\equiv a_na_j\pmod{p}
    \end{equation}
holds.   \hfill $\Box$

\vspace{1mm}

{\it Remark} 17. (cf. \cite[p. 152, Remark 6.1]{al1}). Taking $n=j=0$ in the 
congruence \eqref{con47} gives $a_0\equiv a_0^2(\bmod{\,p})$. This yields 
that either $a_0\equiv 0(\bmod{\,p})$ or $a_0\equiv 1(\bmod{\,p})$.
In the first case, taking $n=0$ and $j\in \{0,1,\ldots, p-1\}$
gives $a_j \equiv 0(\bmod{\,p})$; hence $a_{pn+j}\equiv a_na_j\equiv 0
(\bmod{\,p})$  for all $n$'s and $j$'s. This means that $a_n$ is 
a zero sequence modulo $p$. What precedes implies that such a sequence either
satisfies $a_n=0$ for all $n\ge 0$ or $a_0=1$.\hfill $\Box$
\vspace{1mm}

An analogous definition of double Lucas property is given also 
by McIntosh \cite{m2} as follows: 

\vspace{1mm}

{\it Definition}. The function  $L: \Bbb N\times \Bbb N\rightarrow \Bbb Z$ has 
the {\it double Lucas property} if $L(n,m)=0$ for all $n<m$, 
and for every prime $p$, every 
$n,m\ge 0$, and every $r,s$ with $0\le r,s\le p-1$ the congruence
     \begin{equation}\label{con48}
L(np+r,mp+s)\equiv L(n,m)L(r,s)\pmod{p}
    \end{equation}
holds.\hfill $\Box$
 \vspace{1mm}

Notice that Lucas' theorem (the congruence \eqref{con4}) 
and the congruence \eqref{con52} show that both  functions   
$C(n,m),B(n,m): \Bbb N\times \Bbb N\rightarrow \Bbb Z$ defined as
$C(n,m)={n\choose m}$ and  $B(n,m)={n+m\choose m}$
have the double Lucas property.
McIntosh \cite{m2} presents various properties of the function $L(n,k)$ and 
their connection with tre Lucas property. A typical result
is as follows: {\it If $L(n,m)$ has the double Lucas property, 
then the function $F(n)=\sum_{m=0}^nL(n,m)$ has the Lucas property}.

In 1999 J.-P. Allouche \cite[Proposition 7.1]{al1} proved the following 
result: {\it Let $m$ be a positive integer, let $e_1=2,e_2,\ldots, e_m$ be 
integers such that $e_j\le e_{j+1}\le 2e_j$ for $j=1,2,\ldots,m-1$,
and let $r_1,r_2,\ldots,r_m$ be positive integers. Then the  
sequence $(u_n)_{n\ge 0}$ defined by
  \begin{equation}\label{con49}
u_n={2n\choose n}^{r_1}{e_2n\choose 2n}^{r_2}{e_3n\choose e_2n}^{r_3}
\cdots {e_mn\choose e_{m-1}n}^{r_m}
  \end{equation}
has the Lucas property.}

In particular, {\it if $e_{j+1}-e_j=1$ for all $j=1,2,\ldots,m-1$,
and $r_1,r_2,\ldots ,r_m$ are positive integers,
then the above result implies that the sequence $(u_n)_{n\ge 0}$
defined as 
  $$
u_n={2n\choose n}^{r_1}{3n\choose n}^{r_2}\cdots {(m+1)n\choose n}^{r_m}
  $$
has the Lucas property $($see \cite{m2}$)$}. 

The {\it Ap\'{e}ry numbers} $A_1(n)$ and $A_2(n)$ defined as
      $$
A_1(n)=\sum_{k=0}^n{n\choose k}^2{n+k\choose k}^2,
A_2(n)=\sum_{k=0}^n{n\choose k}^2{n+k\choose k},\,\, n=0,1,\ldots,
    $$
arose in Ap\'{e}ry's proof in 1979 of the  irrationality of  
 $\zeta(3)$ \cite{ap}. $(A_1(n))_{n\ge 0}$ and $(A_1(n))_{n\ge 0}$  are
 Sloane's sequences A005259 and A005258  in \cite{slo}, respectively.

The Ap\'{e}ry numbers modulo a prime were studied in 1982 by I. Gessel  
who proved \cite[Theorem 1]{ge} the following result: 
{\it If $n=n_0+n_1p+\cdots +n_sp^s$ is the $p$-adic expansion of $n$, then 
       \begin{equation}\label{con50}
A_1(n)\equiv \prod_{i=0}^sA_1(n_i)\pmod{p}.
     \end{equation}
In other words, the sequence $(A_1(n))_{n\ge 1}$ has the 
Lucas property.}

Similarly, {\it the sequence $(A_2(n))_{n\ge 0}$ satisfies the 
Lucas property $($see \cite{de}$)$}.

In 2008 Y. Jin, Z-J. Lu and A.L. Schmidt \cite[(ii) of Lemma 2]{jls} 
proved that the sums of powers of binomial coefficients  
have the Lucas property, that is: {\it For a positive integer 
$s$, let $(a_n^{(s)})_{n\ge 0}$ be a sequence defined as
 $$
a_n^{(s)}=\sum_{k=0}^n{n\choose k}^s,\quad n=0,1,2,\ldots.
 $$
Then for every prime $p$, every 
$n\ge 0$, and every $j\in \{0,1,\ldots,p-1\}$ the congruence
     \begin{equation}\label{con51}
a_{pn+j}^{(s)}\equiv a_n^{(s)}a_j^{(s)}\pmod{p}
    \end{equation}
holds.}

\vspace{1mm}
{\it Remark} 18. The above result implies that the residues 
of Pascal's triangle modulo $p$ have a self-similar
structure (see, e.g., \cite{fr2}, \cite[Section 5]{gr} and  
\cite{wo}).\hfill $\Box$
 \vspace{1mm}

For a prime $p$ and a positive integer $k$, in 1994 M. Razpet \cite{raz1} 
considered the $p^k\times p^k$ matrix 
$A(k,p)=[a_{i,j}(k,p)]_{0\le i\le p^k-1}^{0\le j\le p^k-1}$,
whose the entry $a_{i,j}(k,p)$ is defined as the remainder 
of the division of ${i\choose j}$ by $p$. In particular, 
for $k=1$ we write $A(p)=A(1,p)=[a_{i,j}(p)]_{0\le i\le p-1}^{0\le j\le p-1}$.  
M. Razpet \cite{raz1}  noticed that for every $k\ge 1$ 
and every prime $p$, the matrix $A(k,p)$ is the $k$-fold tensor 
(or Kronecker) product
of the matrix $A(p)$ by itself in the field $\Bbb Z_p$,
that is,  
$A(k,p)=\underbrace{A(p)\otimes A(p)\cdots \otimes A(p)}_{k}=A(p)^{\otimes{k}}$.
Note that matrix indices start at index pair $(0,0)$.
This is an algebraic and ``square'' representation of the oft-noted 
self-similarity structure of Pascal's triangle (see, e.g., \cite{hs} 
and \cite{wo}).

Furthermore, as noticed in \cite[p. 378]{raz1}, by Lucas' theorem  we have
       \begin{equation}\label{con54}
a_{i,j}(k,p)\equiv a_{i_0,j_0}(p)a_{i_1,j_1}(p)\cdots 
a_{i_{k-1},j_{k-1}}(p)\pmod{p},
        \end{equation}
{\it where $0\le i,j\le p^k-1$, $i=i_0+i_1p+\cdots +i_{k-1}p^{k-1}$ and 
$j=j_0+j_1p+\cdots +j_{k-1}p^{k-1}$ with $0\le i_l,j_l\le p-1$
for all $l=0,1,\ldots, k-1$.}   

\vspace{1mm}
{\it Remark} 19. In \cite{pr1} M. Prunescu pointed out that Pascal's 
triangle modulo $p^k$ is not a limit of tensor
powers of matrices if $k\ge 2$. However,  Pascal's 
triangle modulo $p^k$ are $p$-automatic, and consequently can be 
produced by matrix substitution and are projections of double sequences
produced by two-dimensional morphisms (see \cite{ash}).  \hfill $\Box$
 \vspace{1mm}

In 2003 D. Berend and N. Kriger \cite[Theorem 5]{bk} 
proved that {\it there exist  uncountably  many infinite 
matrices $A=[a_{i,j}]_{m,n=0}^{\infty}$ satisfying the double Lucas property,
that is the congruences  
      \begin{equation}\label{con55}
a_{m,n}\equiv \prod_{i=0}^ka_{m_i,n_i}\pmod{p}
        \end{equation}
are satisfied for  every prime $p$ and all nonnegative integers $m$ and $n$ 
with  $p$-adic expansions $m=\sum_{i=0}^km_ip^i$ and $n=\sum_{i=0}^kn_ip^i$.}

In 1998 N.J. Calkin \cite{cal} investigated divisibility properties for 
sums of powers of binomial coefficients $f_{n,a}$ defined as
   $$
f_{n,a}=\sum_{k=0}^n{n\choose k}^a,
  $$
where $n$ and $a$ are nonnegative integers. Then $f_{n,0}=n+1$,
$f_{n,1}=2^n$ and $f_{n,2}={2n\choose n}$. 
The  sequences $(f_{n,a})_{n\ge 0}$ for $a=3,4,5,6$ are Sloane's sequences
A000172 (Franel numbers), A005260,
A005261, A069865 in \cite{slo}, respectively. Calkin \cite[Lemma 4]{cal}
proved that {\it for every positive integer  $a$, 
the sequence $(f_{n,a})_{n\ge 0}$ has the Lucas property.
This means that if $p$ is a prime and if 
$n=n_0+n_1p+\cdots +n_sp^s$ is the $p$-adic expansion of $n$, then 
        \begin{equation}\label{con55+}
 f_{n,a}\equiv \prod_{i=0}^s f_{n_i,a}\pmod{p}.
  \end{equation}}
Calkin \cite[p. 21]{cal} also noticed that {\it for any  $a\in\{1,2,\ldots\}$
the  sequence $(h_{n,a})_{n\ge 0}$ defined as
   $$
h_{n,a}=\sum_{k=0}^n(-1)^k{2n\choose k}^a,
  $$
also has the Lucas property.}

 For a positive integer $n$ the central trinomial coefficient
$T_n$ is the largest coefficient in the expansion $(1+x+x^2)^n$
(Sloane's sequence A002426 in \cite{slo}). 
It is easy to express $T_n$ in terms of trinomial coefficients as 
   $$
T_n=\sum_{k\ge 0}{n\choose k,k,n-2k},
   $$ 
where we use the convention that if any multinomial coefficient has a negative
number on the bottom then the coefficient is zero. 
In 2006 E. Deutsch and B.E.  Sagan  \cite{ds} proved  
that  the sequence $(T_n)_{n\ge 0}$ has  the Lucas property. 
Namely,  by 
\cite[Theorem 4.7]{ds} if {\it $p$ is a prime and $n=n_0+n_1p+\cdots +n_sp^s$
is a positive integer with $0\le n_i\le p-1$ for all $i=0,1,\ldots, s$, 
then
      \begin{equation}\label{con56}
T_n\equiv \prod_{i=0}^s T_{n_i}\pmod{p}.
     \end{equation}}
Furthermore, E. Deutsch and B.E.  Sagan  \cite[Theorem 4.4]{ds} proved
the following result for central binomial coefficients ${2n\choose n}$
(Sloane's sequence A000984 in \cite{slo}): 
{\it Let $p$ be a prime and let $n=n_0+n_1p+\cdots +n_sp^s$
be a positive integer with $0\le n_i\le p-1$ for all $i=0,1,\ldots, s$.
For every $j\in\{0,1,\ldots ,p-1\}$ let $\delta_{p,j}(n)$ be the number 
of elements of the set $\{n_0,n_1,\ldots,n_s\}$ equal to $j$. 
Then 
       \begin{equation}\label{con57}
{2n\choose n}\equiv\left\{ \begin{array}{ll}
\prod_{j}{2j\choose j}^{\delta_{p,j}(n)}
&\pmod{p} \quad if\,\, n_i\le p/2\,\ for\,\, all\,\, i=0,1,\ldots s,\\
0   &\pmod{p}\quad otherwise, 
     \end{array}\right.,
\end{equation}
where the summation ranges over all 
$j\in\{0,1,\ldots,p-1\}$ such that $\delta_{p,j}(n)>0$.}

In 2009 M. Chamberland and K. Dilcher \cite{chd1}
studied the divisibility properties of the  sums $u(n)$ defined as
  $$
u(n)=\sum_{k=0}^n(-1)^{k}{n\choose k}{2n\choose k},\quad n=0,1,2,\ldots.
  $$
 Under this  notation, the authors proved \cite[Theorem 2.2]{chd1}
{\it that for every prime $p\ge 3$ and all integers  $m\ge 0$
and $r$ such that $0\le r\le (p-1)/2$ we have}
  \begin{equation}\label{con58}
u(mp+r)\equiv u(m)u(r)\pmod{p}.
  \end{equation}
 As an application, the authors proved \cite[Corollary 2.1]{chd1}  that 
{\it for every prime $p\ge 3$ and every integer 
$n=n_0+n_1p+\cdots +n_sp^s$ with $0\le n_i\le (p-1)/2$ for each 
$i=0,1,\ldots,s$, we have}
  \begin{equation}\label{con59}
u(n)\equiv u(n_0) u(n_1)\cdots u(n_s)\pmod{p}.
  \end{equation}
Similarly, if the sums $w(n)$ are defined as
  $$
w(n)=\sum_{k=0}^{n-1}(-1)^{k}{2n-1\choose k}{n-1\choose k},\quad n=0,1,2,\ldots,
  $$
then by \cite[Corollary 2.2]{chd1}, 
{\it for all primes  $p\ge 3$ and positive integers 
$m$ and $r$ with $(p+1)/2\le r\le p-1$} 
  \begin{equation}\label{con60}
u(mp+r)\equiv w(m+1)u(r)\pmod{p}.
  \end{equation}

  \vspace{1mm}

{\it Remark} 20. We point out that the Lucas property holds for a general 
family  of sequences considered in 2006 by T.D. Noe \cite{noe}. \hfill$\Box$ 
   \vspace{1mm}

For all nonnegative integers $i$ and $j$ let $w(i,j|a,b,c)$ denote the number 
of all paths in the plane from $(0,0)$ to $(i,j)$ with steps $(1,0)$, $(0,1)$,
$(1,1)$, and with positive integer weights $a$, $b$, $c$, respectively. 
The explicit formula for $w(i,j|a,b,c)$ was obtained by several
authors by using combinatorial arguments (see, e.g., \cite{fr}):
   $$
w(i,j|a,b,c)=\sum_{k\ge 0}{k\choose i}{i\choose k-j}a^{k-j}b^{k-i}c^{i+j-k}.
  $$
Actually, $k$ in the above sum runs from $\max \{i,j\}$ to $i+j$.
In the case $a=b=c=1$, we have even  the {\it Delannoy numbers} which count
the usual, unweighted lattice paths from the point $(0,0)$ to the point
$(i,j)$ with steps along the vectors $(1,0)$, $(0,1)$ and $(1,1)$.
If $i=j=n$, then the numbers $w(n,n|1,1,1)$, $n=0,1,2,\ldots$ are called the
{\it central Delannoy numbers} (Sloane's sequence A001850  in \cite{slo}).  

 In 2002  M. Razpet \cite[Theorem 1]{raz2} proved the following double Lucas 
property of $w(i,j|a,b,c)$:
{\it Let $p$ be a prime and let $\alpha,\beta,\gamma,\delta$ be nonnegative 
integers where $0\le \beta <p$ and $0\le \delta<p$. Then the congruence 
   \begin{equation}\label{con61} 
w(\alpha p+\beta, \gamma p+\delta|a,b,c)\equiv 
w(\alpha,\gamma |a,b,c)w(\beta,\delta |a,b,c)\pmod{p}
  \end{equation}
holds for all positive integers $a,b,c$}.

  \vspace{1mm}
{\it Remark} 21. Razpet \cite{raz2} notice that the congruence 
\eqref{con61} is particularly  true for the {\it Delannoy numbers} 
$D(i,j):=w(i,j|1,1,1)$ as proven in another way in 1990 by M. Razpet 
\cite{raz3} and by M. Sved and R.J. Clarke \cite{scl} 
(see also \cite{ds} and \cite{ely}). \hfill$\Box$ 
  \vspace{1mm}

In 2004 H. Pan \cite[Theorem 1]{pa1} proved the following result:
{\it Suppose $\lambda(x_1,\ldots,x_n)=\sum_{\Phi\not= 
I\subseteq \{1,\ldots,n\}}\alpha_I\prod_{i\in I}X_i$,
$\alpha_I\in\Bbb F$, is a polynomial over the 
finite field  $\Bbb F$ with $q$ elements. Let 
$w_{\lambda}(k_1,\ldots,k_n)$ be the coefficient of 
$\prod_{i=1}^nX_i^{k_i}$ in the formal power series 
$\frac{1}{1-\lambda(x_1,\ldots,x_n)}$. Then $w_{\lambda}$ satisfies the 
double  Lucas property, i.e., for any nonnegative integers $a_1,\ldots,a_n$
and $0\le b_1,\ldots,b_n<q$, 
  \begin{equation}\label{con62}
w_{\lambda}(a_1q+b_1,\ldots,a_nq+b_n)=w_{\lambda}(a_1,\ldots,a_n)
w_{\lambda}(b_1,\ldots,b_n).
  \end{equation}}

  \vspace{1mm}
{\it Remark} 22. If $p$ is a prime and  $\Bbb F$ is the field 
$\Bbb Z_p=\{0,1,\ldots,p-1\}$, then the equality ``$=$'' in  
\eqref{con62} becomes $\equiv 0(\bmod{\,p})$.    \hfill$\Box$ 

\subsection{Further Lucas type congruences}\label{subsec4.2}

For nonnegative integers $n$ and $k$  
{\it Stirling numbers of the second kind} ${n\brace k}$
(Sloane's sequence A008277 in \cite{slo}) are recursively defined as:
   $$
{n\brace k}=\left\{\begin{array}{ll}
1 & {\rm if}\,\, n=0, k=0,\\
0 & {\rm if}\,\, n>0, k=0,\\
0 & {\rm if}\,\, n=0, k>0,\\
k{n-1\brace k}+{n-1\brace k-1} & {\rm if}\,\, n>0, k>0.
  \end{array}\right.
   $$
${n\brace k}$ presents the number of ways of partitioning a set of $n$ elements 
into $k$ nonempty sets (i.e., $k$ set blocks). 
They (as well as Stirling numbers of the first kind defined below) 
are named after James Stirling, who introduced them in 1730 \cite{sti}. 

In 1988 M. Sved \cite[p. 61, Theorem]{sv} showed the following result:
{\it Let $n$ and $m$ be  nonnegative integers, and let $p$ be a an odd prime
such that $p$ does not divide $m$. Put
  $$
n'=\left\lfloor\frac{pn-p\lfloor m/p\rfloor -1}{p-1}\right\rfloor,
  $$
and let $n'=\sum_{i=0}^hn_i'p^i$ and $m=\sum_{i=0}^hm_ip^i$
be the expansions of $n'$ and $m$ to base $p$. Then 
    \begin{equation}\label{con63}
{n\brace m}\equiv {n_0'\brace m_0}\prod_{j=1}^h{n_j'\choose m_j}.      
  \end{equation}}

In 2000 R. S\'{a}nchez-Peregrino \cite[Proposition 3.1]{pe}  
proved that  {\it if $n,m,r$ and $s$ are nonnegative integers such that 
$0\le s\le r\le p-1$ and $m\le n\le p-1$, then} 
   \begin{equation}\label{con64}
{np+r\brace mp+s}\equiv {n-m+r\brace s}{n\choose m}+
{n-m+r+1\brace s+p}{n\choose m-1}\pmod{p}.
  \end{equation}

{\it Notice also that under the hypothesis that $r+n-m+1<s+p$, the  
congruence \eqref{con64} reduces to} 
 \begin{equation}\label{con65}
{np+r\brace mp+s}\equiv {n-m+r\brace s}{n\choose m}\pmod{p}.
  \end{equation}

Furthermore, by \cite[Proposition 4.1]{pe}, {\it if $r,s,a$ and $f$ 
are nonnegative inegers, then}
 \begin{equation}\label{con66}
{ap^f+r\brace s}\equiv \sum_{i_0+i_1+\cdots +i_f}{a\choose 
i_0,i_1,\ldots ,i_f}{r+i_0\brace s-\sum_{l=1}^fi_lp^f}\pmod{p}.
 \end{equation}
  \vspace{1mm}

{\it Remark} 23. As noticed in \cite[Remark 3.1]{pe}, in the case 
$r<p$ the congruence \eqref{con64} gives the formulas (4.17) and (4.18) of 
F.T. Howard \cite{ho} from 1990.\hfill$\Box$
   \vspace{1mm}

For nonnegative integers $n$ and $k$  
{\it Stirling numbers of the first kind} ${n\brack k}$
(Sloane's sequence A008275 in \cite{slo}) are defined by the recurrence 
relation 
 $$
{n\brack k}=\left\{\begin{array}{ll}
1 & {\rm if}\,\, n=0, k=0,\\
0 & {\rm if}\,\, n>0, k=0,\\
0 & {\rm if}\,\, n=0, k>0,\\
(n-1){n-1\brack k}+{n-1\brack k-1} & {\rm if}\,\, n>0, k>0.
  \end{array}\right.
   $$
The absolute value  of ${n\brack k}$ 
(Sloane's sequence A094638 in \cite{slo}) denotes, as usual, the number of 
permutations of $n$ elements which contain exactly $k$ permutation cycles. 
 
In 1993  R. Peele, A.J. Radcliffe and H.S. Wilf \cite[Proposition 2.1]{prw}
proved the following analogue of Lucas' theorem for the numbers ${n\brack k}$:
{ \it Let $p$ be a prime and let $n$ and $k$ be integers with $1\le k\le n$.
Let $n'=\lfloor n/p\rfloor$ and $n_0=n-n'p$. Further, define integers $i$ and 
$j$ as follows:
  $$
k-n'=j(p-1)+i\quad with\,\, 0\le i<p-1\,\, if\,\, n_0=0;
\,\, 0< i\le p-1\,\, if\,\, n_0>0.
  $$
Then 
  \begin{equation}\label{con67}
{n\brack k}\equiv (-1)^{n'-j}{n'\choose j}{n_0\brack i}\pmod{p}.
  \end{equation}}

For a nonnegative integer  $k$ let $J_k(z)$ be the {|it Bessel function of the 
first kind}. Put
     $$
f_k(z)=\frac{J_k(2\sqrt{z})}{z^{k/2}}=
\sum_{i=0}^{\infty}\frac{(-1)^iz^i}{i!(i+k)!}.
    $$
Furthermore, define the polynomial $u_i(k;x)$ by means of
   $$
\frac{k!f_k(xz)}{f_k(z)}=\sum_{i=0}^{\infty}u_i(k;x)
\frac{z^i}{i!(i+k)!}.
   $$   
Certain Lucas type congruences for $w_i(x)=u_i(0;x)$ and the 
integers $w_i=w_i(0)$ with $i=0,1,2\ldots,$ were derived by 
L. Carlitz \cite{ca2} in 1955, and an interesting application was presented 
($(w_n)_{n\ge 0}$ is Sloane's sequence A000275). 
   In 1987 F.T. Howard \cite[Theorem 1]{ho4} proved a more general result as 
follows: {\it Let  $k$, $n$ and $s$ be nonnegative integers, 
and let $p$ be a prime such that $p\ge 2k$ and  $s<p-2k$. Then the 
numbers $u_i(k):=u_i(k;0)$ are integral $(\bmod{\,p})$ for 
all $i=0,1,2\ldots;$ in particular, $u_n(0)$ and $u_n(1)$ are 
positive integers for  all $n=0,1,2\ldots$. Furthermore, 
for any fixed $k\ge 0$ and every prime $p$   the congruence
   \begin{equation}\label{con68}
  u_{np+s}(k)\equiv u_s(k)\cdot w_n\pmod{p}
  \end{equation}
holds for all $n\ge 0$ and $0\le s\le n-1$}.

{\it With the assumptions of the above statement, if $m$ is a nonnegative 
integer with the expansion $m=\sum_{i=0}^sm_ip^i$ to base  $p$
satisfying $m_0<p-2k$, then the congruence \eqref{con69} with $k=0$ implies 
Carlitz's result  \cite{ca2} from 1955 which asserts that the sequence 
$(w_n)_{n\ge 0}$ has the Lucas property, i.e., 
    \begin{equation}\label{con69}
w_m\equiv \prod_{i=0}^sw_{m_i}\pmod{p},
   \end{equation}}
Furthermore, {\it the following  two  congruences are satisfied
\cite[p. 306, Corollary and Theorem 2]{ho4}:
      \begin{equation}\label{con70}
u_m(k)\equiv u_{m_0}(k) \prod_{i=1}^sw_{m_i}(0)\pmod{p},
   \end{equation}
and
     \begin{equation}\label{con71}
u_{np-k}(k)\equiv (-1)^ku_0(k)\cdot w_n(0)\pmod{p}.
   \end{equation}}

Let $p$ be a prime and let $n,r,l$ and $a$ be positive integers.
Following Z.-W. Sun and D. Wan \cite{sw},  the {\it normalized 
cyclotomic $\psi$-coefficient} is defined as
  \begin{equation}\label{con72}
{n\brace r}_{l,p^a}:=p^{-\left\{\lfloor \frac{n-p^{a-1}-lp^a}{\varphi(p^a)}
\right\{\rfloor}\sum_{k\equiv r\,(\bmod{\,p^a})}(-1)^k{n\choose k}
{(k-r)/p^a\choose l}.
  \end{equation}
In 2008 Z.-W. Sun and D. Wan  \cite[Theorem 1.1]{sw} proved that 
{\it if $p$ is any prime, $r$ is an integer
and $a,l,n,s,t$ are positive integers with $a\ge 2$ and $s,t<p$, then}
  \begin{equation}\label{con73}
{pn+s\brace pr+t}_{l,p^{a+1}}\equiv (-1)^t{s \choose t}
{n\brace r}_{l,p^a}\pmod{p}.
    \end{equation}

It is noticed in \cite[Remark 1.1]{sw} that in the case $l=0$ the congruence
\eqref{con73} is equivalent to Theorem 1.7 in \cite{sd} due to Z.-W. Sun and 
D.M. Davis in 2007. Under the same conditions
preceeding the congruence \eqref{con73}, Sun and Davis 
\cite[Theorem 1.7]{sd} proved the following congruence of Lucas' type: 
   \begin{equation}\begin{split}\label{con74}
&\frac{1}{\lfloor n/p^{a-1}\rfloor !}\sum_{k\equiv r\,(\bmod{\,p^a})}(-1)^{pk}
{pn+s\choose pk+t}\left(\frac{k-r}{p^{a-1}}\right)^l\\
\equiv & \frac{1}{\lfloor n/p^{a-1}\rfloor !}\sum_{k\equiv r\,(\bmod{\,p^a})}
(-1)^{k}{n\choose k}{s\choose t}\left(\frac{k-r}{p^{a-1}}\right)^l\pmod{p}.
     \end{split}\end{equation}

J. Boulanger and J.-L. Chabert \cite{bc} have extended Lucas' theorem 
to Linear Algebra and Even Topology. Their  result can be 
briefly exposed as follows.  Let $V$ be a discrete valuation 
domain with finite residue field. Denote by
$K$ the quotient field of $V$, by $v$ the corresponding valuation 
of $K$, by $\frak{m}$ the maximal ideal of $V$, and by $q$ 
the cardinality of the residue field $V/\frak{m}$. We denote 
by $\widehat{K}$, $\widehat{V}$ and $\widehat{\frak{m}}$ the completions
of $K$, $V$ and $\frak{m}$, respectively,   with respect
to the $\frak{m}$-adic topology and we still denote by $v$ the extension 
of $v$ to $\widehat{K}$. Consider the ring ${\rm Int}(V)$ of integer-valued 
polynomials on $V$, that is,
  $$
{\rm Int}(V)=\{ f\in K[X]:\, f(V)\subseteq V\}.
  $$
A basis $C_n(X)$ of the $V$-module ${\rm Int}(V)$ can be constructed as follows 
\cite[Chapter II, \S 2 ]{cc}. We choose a generator $t$ of $\frak{m}$
and a set $U=\{u_0=0,u_1,\ldots, u_{q-1}\}$ of representatives of $V$
modulo $\frak{m}$. It is known that each element $x$ of $\widehat{V}$
has a unique $t$-adic expansion 
$$
x=\sum_{j=0}^{\infty}x_jt^j\quad {\rm with}\,\, x_j\in U\,\, 
{\rm for\,\, each} \,\, j\in \Bbb N.
 $$

We now construct a sequence $(u_n)_{n\ge 0}$ of elements of $V$ which will
replace the sequence of nonnegative integers. Taking $q$ as the basis of the
numeration, that is, writing every positive integer $n$ in the form
$n=\sum_{i=0}^kn_iq^i$ with $0\le n_i<q$ for each $i=0,1,\ldots,k$,
we extend the sequence $(u_j)_{0\le j<k}$ in the following way:
  $$
u_n=u_{n_0}+u_{n_1}t+u_{n_2}t^2+\cdots +u_{n_k}t^k.
 $$
We then replace the binomial polynomials 
  $$
{X\choose n}=\frac{X(X-1)(X-2)\cdots (X-n+1))}{n!}
 $$
$($which form a basis of the $\Bbb Z$-module ${\rm Int}(\Bbb Z)
=\{f\in \Bbb Q[X]:\, f(\Bbb Z)\subseteq \Bbb Z\}$ of integer-valued 
polynomials on $\Bbb Z$$)$ by the polynomials defined as
  $$
C_n(X)=\prod_{k=0}^{n-1}\frac{X-u_k}{u_n-u_k},\,\, n=1,2,\ldots,\,\, {\rm and}
\,\, C_0=1.
  $$ 
Then by \cite[Theorem II.2.7]{cc}, the sequence of  polynomials $(C_n(X))_{n\ge 0}$ 
form a basis of the $V$-module ${\rm Int}(V)$. 
In 2001 J. Boulanger and J.-L. Chabert \cite[Theorem 2.2]{bc}
proved the following ``generalized Lucas' theorem'':
{\it If 
$$
n=n_0+n_1q+\ldots +n_kq^k
$$ 
is the $q$-adic expansion of a positive
integer $n$, and if 
$$
x=x_0+x_1t+\ldots +x_jt^j+\ldots
$$ 
is the $t$-adic expansion of an element $x$ of $\widehat{V}$, then} 
     \begin{equation}\label{con75}
C_n(x)\equiv C_{n_0}(x_0)C_{n_1}(x_1)\cdots C_{n_k}(x_k) 
\pmod{\widehat{\frak{m}}}. 
  \end{equation}

  \vspace{1mm}
{\it Remark} 24.
Notice also that in  1993 N. Zaheer \cite{za} generalized Lucas' theorem to 
vector-valued abstract polynomials in vector spaces.  \hfill$\Box$ 
  
\section{Lucas type theorems for some generalized binomial 
coefficients}\label{sec5}

\subsection{Generalized binomial coefficients and related Lucas type 
congruences}\label{subsec5.1}

Let $A$ and $B$ be nonzero integers. The {\it Lucas sequence} 
$u_0,u_1,u_2,\ldots$ is defined recursively as 
       \begin{equation}\label{con76}
u_0=0,u_1=1\quad {\rm and} \quad u_{n+1}=Au_n-Bu_{n-1}\quad {\rm for}\,\,
n=1,2,3,\ldots .
    \end{equation} 
The companion sequence of Lucas sequence $(u_n)_{n\ge 0}$
is the sequence $(v_n)_{n\ge 0}$ recursively defined as 
       \begin{equation}\label{con77}
v_0=2,v_1=A\quad {\rm and} \quad v_{n+1}=Av_n-Bv_{n-1}\quad {\rm for}\,\,
n=1,2,3,\ldots .
    \end{equation} 
It is well known that for all $n=0,1,2,\ldots$
  $$
u_n=\frac{\alpha^n-\beta^n}{\alpha -\beta}\quad {\rm and}\quad
v_n=\alpha^n+\beta^n,
  $$
where 
  $$
\alpha=\frac{A+\sqrt{\Delta}}{2}, \beta=\frac{A-\sqrt{\Delta}}{2}
\quad {\rm and}\quad \Delta=A^2-4B.
  $$
In fact, $\alpha$  and $\beta$ are roots of the characteristic equation 
$x^2-Ax+B=0$.   
Note that for $A=1,B=-1$ the terms of the sequence
$(u_n)_{n\ge 0}$ defined by  \eqref{con76} are the well-known 
{\it Fibonacci numbers} $F_n$ defined recursively as $F_0=0$, $F_1=1$  and 
 $$
F_{n+1}=F_n+F_{n-1}\quad {\rm for}\,\, n\ge 1.
 $$ 
Fibonacci numbers are in fact the  Lucas sequence $(u_n)_{n\ge 0}$ given by 
\eqref{con76} with $u_0=0$ and  $u_1=1$.

Similarly, the {\it Lucas numbers} $L_n$ are defined by $L_0=2$,  $L_1=1$ and 
$$
L_{n+1}=L_n+L_{n-1}\quad {\rm for}\,\, n\ge 1.
$$ 
Fibonacci numbers $F_n$ and Lucas numbers $L_n$ are  given as Sloane's 
sequences A000045 and  A000032 in \cite{slo}, respectively.

Let $a:=(a_n)_{n\ge 0}$ be a sequence of real or complex numbers such that 
$a_n\not= 0$ for all $n\ge 1$. The $a$-{\it nomial coefficients}
  (or the {\it  generalized binomial coefficients})
(associated to the sequence $a$)  are defined by 
  $$
{n\brack k}_a=\frac{a_na_{n-1}\cdots a_1}{(a_ka_{k-1}\cdots a_1)(a_{n-k}
a_{n-k-1}\cdots a_1)}
\quad{\rm for}\,\, n\ge 2\,\, {\rm and}\,\, 1\le k\le n-1,
  $$
and 
  $$
{n\brack 0}_a={n\brack n}_a=1\quad{\rm for}\,\, n\ge 0.
   $$
This definition was suggested  in  1915 by Georges Fonten\'{e} 
in his one-page note \cite{fon}. 
A number of authors have considered different classes of  generalized 
binomial coefficients ${n\brack k}_a$ 
(usually, when  $a:=(a_n)_{n\ge 0}$ is an integer sequence). Related investigations 
were done in 1913 by R.D. Carmichael \cite{car}, in 1936 by  M. Ward, \cite{wa}, 
in 1967 by R.D. Fray \cite{fr2} and V.E. Hoggatt \cite{hog}, in 1969 by H.W. 
Gould \cite{go}, and later by several authors 
(\cite{hol1}, \cite{hol3}, \cite{kw1}, \cite{kw}, \cite{noe}, \cite{tf} and 
\cite{to}). For example, in 1989 D.E. Knuth and H.S. Wilf 
\cite[Proposition 3]{kw} generalized Kummer's  theorem for the $a$-nomial 
coefficients ${m+k\brack m}_a$, where $a=(a_n)_{n\ge 1}$ is a 
sequence of positive integers. Consequently, they obtained 
\cite[Theorems 1 and 2]{kw} Kummer's  theorem for the Gaussian $q$-nomial 
coefficients ${m+k\brack m}_q$ where  $q>1$ is an integer and 
for the  Fibonomial coefficients ${m+k\brack m}_{\mathcal F}$ defined below, 
respectively. 

In general, even if the all terms of a sequence $a=(a_n)_{n\ge 0}$  are 
integers,  ${n\brack k}_a$ may not be integers. 
In 1913 R.D.  Carmichael \cite[page 40]{car} proved   that {\it if the 
sequence  $a:=(a_n)_{n\ge 1}$  of positive integers is defined  recursively as 
  $$
a_1=a_2=1,\quad{\rm and}\quad a_{n+1}=ca_n+da_{n-1}\quad {\rm for}\,\, 
n=2,3,4,\ldots,
  $$
where $c$ and $d$ are integers, then the all $a$-nomial coefficients 
are integers}. For a more general result see Remark 28.

If $u:=(u_n)_{n\ge 0}$ is the Lucas sequence defined by \eqref{con76}, 
and if $A\not= \pm 1$ or $B\not= 1$, then $u_1,u_2,\ldots$
are nonzero (see, e.g.,  \cite{hsu}), and so are $v_1=u_2/u_1$, 
$v_2=u_4/u_2,\ldots$,  where $v:=(v_n)_{n\ge 0}$ is the companion sequence 
of the sequence $(u_n)_{n\ge 0}$ given by \eqref{con77}. 
In the case when $A^2=B=1$, 
then as noticed in \cite{hsu} $u_n=0$ if and only if $3\mid n$.
If $v_n=0$, then $u_{2n}=u_nv_n=0$; hence $3\mid n$ and $u_n=0$,
which is impossible since $v_n^2-\Delta u_n^2=4B^n$ (cf. \cite{hu}).
Thus $v_0,v_1,v_2,\ldots$ are all nonzero.

If f $A\not= \pm 1$ or $B\not= 1$ the {\it Lucas $u$-nomial coefficient} 
${n\brack k}_u$ with  $1\le k\le n$ is the generalized 
binomial coefficient associated to the Lucas sequence  $u:=(u_n)_{n\ge 0}$
defined by \eqref{con76}, that is,
   $$
{n\brack k}_u=\frac{u_nu_{n-1}\cdots u_1}{(u_ku_{k-1}\cdots
u_1)(u_{n-k}u_{n-k-1}\cdots u_{1})}\quad{\rm for}\,\,n\ge 2\,\,{\rm and}\,\,
1\le k\le n-1, 
  $$
and ${n\brack 0}_u={n\brack n}_u=1$ for all $n\ge 0$. 

In the sam way  we define the $v$-{\it nomial generalized binomial 
coefficient} ${n\brack k}_v$, where $v:=(v_n)_{n\ge 0}$ is the companion 
sequence of the Lucas sequence $(u_n)_{n\ge 0}$ defined by \eqref{con77}.

   \vspace{1mm}
{\it Remark} 25. In the case $A=2$ and $B=1$, \eqref{con76} yields
 $u_n=n$ for all $n=0,1,2,\ldots$, and hence ${n\brack k}_{u}$ is exactly 
the binomial coefficient ${n\choose k}$.\hfill $\Box$
 \vspace{1mm}

Similarly, the {\it Fibonomial coefficients} (or  {\it Fibonacci 
coefficients})  are defined as the generalized binomial coefficients 
associated to the sequence $(F_n)_{n\ge}$ of Fibonacci numbers,  that is,
  $$
{n\brack k}_{\mathcal F}=\frac{F_nF_{n-1}\cdots F_1}{(F_kF_{k-1}\cdots
F_1)(F_{n-k}F_{n-k-1}\cdots F_{1})}\quad{\rm for}\,\,n\ge 2\,\,{\rm and}\,\,
1\le k\le n-1,
  $$
and ${n\brack 0}_{\mathcal F}={n\brack n}_{\mathcal F}=1$ for all $n\ge 0$. 

The Fibonomial coefficients  and the Lucas $u$-nomial coefficients  
were introduced in 1878 by \'{E}. Lucas \cite[\S 9]{l2}, 
and later they have been studied by several authors (see \cite{go}, 
\cite{hog}, \cite{hol2}, \cite{we1}, \cite{hsu} and \cite{hu}). 

The triangle of Fibonomial coefficients is given as  Sloane's sequence 
A010048 in \cite{slo}. It is known (see, e.g., \cite[the equality (D), page 
386]{hog}) that 
  $$
{n\brack k}_{\mathcal F}=F_{k+1}{n-1\brack k}_{\mathcal F}+
F_{n-k-1}{n-1\brack k-1}_{\mathcal F},\,\, {\rm for\,\,}  
0\le k\le n-1, 
  $$
whence by induction immediately follows that  the all
Fibonomial coefficients are integers.

When $A=q+1$ and $B=q$ related to  the sequence 
defined by \eqref{con76}, where $q$ is an integer such that $|q|>1$, 
${n\brack k}_u$, then it coincides with the {\it Gaussian $q$-nomial 
coefficient} ${n\brack k}_q$ because $u_j=(q^j-1)/(q-1)$ for 
$j=1,2,\ldots$, and hence,
  $$
{n\brack k}_q=\frac{(q^n-1)(q^{n-1}-1)\cdots (q^{n-k+1}-1)}
{(q^k-1)(q^{k-1}-1)\cdots (q-1)}.
  $$
The numbers ${n\brack k}_q$ were introduced in 1808 by Gauss 
\cite[\S 5]{ga2}.
It is well known that these numbers satisfy the recursion formula
  $$
{n\brack k}_q=q^k{n-1\brack k}_q+{n-1\brack k-1}_q, \,\, {\rm for}\,\,
0\le k\le n-1.
  $$
The triangles of Gaussian $q$-nomial coefficients for $q=-2,2,3,4,5,6,7,8,9$  
are given as  Sloane's sequences A015109,  A022166, A022167, A022168, A022169,
A022170, A022171, A022172 and A022173 in \cite{slo}, respectively.

 It is easy to see that if $0\le m\le n$, then
 $$
\lim_{q\to 1}{n\brack m}_q={n\choose m}, 
  $$
$$
{n\brack m}_q={n\brack n-m}_q\quad ({\rm symmetry})
 $$
and 
$$
{n\brack m}_q={n-1\brack m-1}_q+q^m{n-1\brack m}_q,
  $$
whence easily follows by induction that if $q$ is any positive integer,
then ${n\brack m}_q$ are also integers for all $n$ and $m$.\hfill$\Box$

 \vspace{1mm}
{\it Remark} 26.  An analogy to the Lucas $u$-nomial coefficients   
${n\brack k}_u$ was obtained in 1995  by W.A. Kimball and W.A. Webb
\cite{kw1} and in 1998 by B. Wilson \cite{wi} in some special cases,
and in 2001 by H. Hu and Z.-W. Sun  \cite{hsu} for the general case 
(see Subsection \ref{subsec5.2}). 
\hfill $\Box$
 \vspace{1mm}

It is known (see, e.g., \cite{lw}, \cite{we1}) that 
the generalized base for the Fibonacci sequence is 
  $$
\mathcal P=\{r_0,r_1,r_2,r_3,r_4,\ldots ,\}=
\{1,3,6,6,12,\ldots,\}
 $$ 
in the sense that any positive  integer $n$ can be uniquely expressed
as 
   $$
n=(n_sn_{s-1}\ldots n_1n_0)_{\mathcal P}:=
n_0+n_1r_1+\cdots +n_{s-1}r_{s-1}+n_sr_s,
  $$
where $0\le n_i< r_{i+1}/r_i$ for each $i=0,1,\ldots,s-1$.

Under the above notations, in 1994 D.L. Wells \cite[Theorem 2]{we1}  
proved that 
 \begin{equation}\label{con78}
{n\brack k}_{\mathcal F}\equiv 
{n_0\brack k_0}_{\mathcal F}\cdot\prod_{i\ge 1}{n_i\choose k_i}\pmod{2}.
  \end{equation}

In 1988 M. Sved \cite{sv} establihed that the geometry of the 
binomial arrays of Pascal's triangle modulo $p$ gives 
a simple interpretation of Lucas' theorem. 
Moreover, as noticed in \cite[p. 58]{sv}, this interpretation can be extended 
to arrays of other combinatorial functions; in particular, Lucas' theorem can 
be generalized to the Gaussian $q$-nomial coefficients as follows.
 {\it Let $p$ be a prime, $q>1$  a positive integer not divisible by $p$, 
and let $a\not= 1$ be the minimal exponent for which
$q^a\equiv 1\,(\bmod{\,p})$; then by Fermat little theorem it follows 
that $a\mid (p-1)$. Further, if  $n=Na+n_0$, $m=Ma+m_0$ with 
$0\le n_0<a$ and $0\le m_0<a$, then \cite[p. 60]{sv}
   \begin{equation}\label{con79}
{n\brack m}_q\equiv {N\choose M}{n_0\brack m_0}_q\pmod{p}.
\end{equation}}

 \vspace{1mm}

{\it Remark} 27. In the same area of research A. B\`{e}s \cite{be}
generalized Lucas' theorem. This accomplishment obviously serves to improve
the security of cryptographic applications modulo prime powers 
\cite{be}. \hfill$\Box$
 \vspace{1mm}

{\it Definition}.  For a positive integer $d$, the  {\it rank of apparition} 
$r=r(d)$ with respect to the integer  sequence $(a_n)_{n\ge 0}$ is 
the least index $n$ for which $d$ divides $a_n$, that is, 
$r(d)=\min \{n\in \Bbb N:\, d\mid a_n\}$ (if $d$ does not divide any $a_n$,
then $r(d)=\infty$).\hfill $\Box$

 \vspace{1mm}
{\it Remark} 28. Let $a=(a_n)_{n\ge 0}$ be an integer sequence. In order 
to guarantee that the all $a$-nomial coefficients ${n\brack k}_a=0$
are integers, it is usually required that the sequence $a=(a_n)_{n\ge 0}$
be {\it regularly divisible}, that is, $p^i\mid a_j$ if and only if 
$r(p^i)\mid j$ for all $i\ge 1$, $j\ge 1$, and all primes $p$. Here 
$r(p^i)$ denotes the rank of  apparition og $p^i$ as defined above. 
The principal class of sequences which are known to be regularly
divisible are the Lucas sequences given by \eqref{con76} for which
$\gcd (A,B)=1$ (see \cite{hsk}).  \hfill$\Box$
 \vspace{1mm}

In 2000  J.M. Holte \cite[Theorem 1]{hol1} proved the following result:
{\it Let $p$ be a prime and let $m$ and $n$ be nonnegative integers. 
Let $r$ be  the rank of apparition of $p$ with respect to the 
Lucas sequence $u=(u_n)$, let $\tau$ be the period of $(u_n)$ modulo
$p$, and let $t=\tau/r$ $($$t$ is necessarily a positive integer$)$.
Furthermore, for $i,j\ge 0$ and for $0\le k,l<r$, let $A_{i,j}(k,l)$
denote the solution of the modulo $p$ recurrence relation
  $$
A_{i,j}(k,l)\equiv   u_{ir+k+1}A_{i,j}(k,l-1)+bu_{jr+l-1}A_{i,j}(k-1,l)
\pmod{p},
  $$
and let $H_{i,j}(k,l)=u_{r+1}^{rij}A_{i,j}(k,l)$. 
Set $n_0=n(\bmod{\, r})$,  $m_0=m(\bmod{\, r})$, $n'=n+r$, $m'=m+r$, 
 $n''=n'(\bmod{\, t})$, and $m''=m'(\bmod{\, t})$. Then
  \begin{equation}\label{con80}
{m+n\brack n}_{u}\equiv {m'+n'\choose n'}H_{m'',n''}(m_0,n_0)\pmod{p}.
  \end{equation}} 
 
Using the above result, with the same notations as above,
Holte \cite[Theorem 3]{hol1} also proved the following result:
{\it Let  $(u_n)$ be the Lucas sequence defined by \eqref{con76},  let 
$p$ be a prime such that $B$ is not divisible by $p$. Set
$\lambda=\max\{0,m''+n''-(p-1) \}$, $n^*=n(\bmod{\, t})$ and 
$m^*=m(\bmod{\, t})$. Then 
   \begin{equation}\label{con81}
{m+n\brack n}_{u}\equiv 
{m'+n'\choose n'}{m''+n''+\lambda t,\choose n''+\lambda t}^{-1}
{m^*\brack n^*+\lambda \tau}_{u}\pmod{p}.
  \end{equation}
  Thus, except when $s=p-1$ and $m''+n''\ge p$, then 
   \begin{equation}\label{con82}
{m+n\brack n}_{u}\equiv 
{m'+n'\choose n'}{m''+n''+\lambda t,\choose n''}^{-1}
{m^*\brack n^*}_{u}\pmod{p}.
  \end{equation}}

 Holte \cite[Section 7]{hol1} noticed that by means of 
a bit of translation, the congruence \eqref{con81} may be transformed into
the following  result obtained in 1992 by D. Wells  \cite{we3} 
(also see \cite{we4}): {\it Let $N=n+m$, and correspondingly, 
$N_0=N(\bmod{\, r})$, $N'= \lfloor N/r\rfloor$, and 
$N''=N'(\bmod{\, s})$. Let $N'=\sum_{j=0}^lN_jp^j$ and 
$m'=\sum_{j=0}^lm_jp^j$ be the $p$-adic expansions of $N'$ and $m'$.
If $p$ is a prime such that $B$ is not divisible by $p$, then 
under the same definitions of $B$ and $t$ as above,  for $N''\ge m''$,
    \begin{equation}\label{con82+}
{N\brack m}_{u}\equiv {N''\choose m''}^{-1}\prod_{j=0}^l{N_j\choose m_j}
{Nr+N_0\brack m''r+m_0}_{u}\pmod{p},
   \end{equation}
and for $N''<m''$,
   \begin{equation}\label{con82++}
{N\brack m}_{u}\equiv \left\{ \begin{array}{ll}
{s+N''\choose m''}^{-1}\prod_{j=0}^l{N_j\choose m_j}
{t+N''r+N_0\brack m''r+m_0}_{u}\pmod{p} & if\,\, s<p-1\\
{s\choose m''}^{-1}\prod_{j=0}^l{N_j\choose m_j}
{(N''+1)t+N''r+N_0\brack m''r+m_0}_{u}\pmod{p} & if\,\, s=p-1.
   \end{array}\right.
    \end{equation}}

\vspace{1mm}
{\it Remark} 29. In 2002  E.R. Tou \cite[Theorem 4]{to} generalized the 
congruence \eqref{con81} modulo product of a finite number of distinct primes. \hfill$\Box$ 

\subsection{Lucas type congruences for some classes of 
Lucas $u$-nomial coefficients   }\label{subsec5.2}

In 2001  H. Hu and Z.-W. Sun \cite[Theorem]{hsu} proved the following result
for the Lucas $u$-nomial coefficients: {\it Let $u=(u_n)_{n\ge 0}$ be 
a Lucas sequence defined by \eqref{con76}. Suppose that $\gcd(A,B)=1$, 
and $A\not= \pm 1$ or $B\not= \pm 1$. Then $u_k\not= 0$ for every $k\ge 1$. 
Let $q$ be a positive integer, let $m$ and $n$ be nonnegative 
integers, and let $R(q)=\{0,1,\ldots, q-1\}$. If $s,t\in R(q)$ then
    \begin{equation}\label{con83}
{mq+s\brack nq+t}_u\equiv {m\choose n} {s\brack t}_u
u_{q+1}^{(nq+t)(m-n)+n(s-t)}\pmod{w_q},
   \end{equation}
where $w_q$ is the largest divisor of $u_q$ relatively prime to 
$u_1,\ldots,u_{q-1}$. If $q$ or $m(n+t)+n(s+1)$ is even, then
    \begin{equation}\label{con84}
{mq+s\brack nq+t}_u\equiv {m\choose n} {s\brack t}_u
(-1)^{(mt-ns)(q-1)}B^{\frac{q}{2}((nq+t)(m-n)+n(s-t))}\pmod{w_q}.
   \end{equation}}

 \vspace{1mm}
{\it Remark} 30. (\cite[Remark 1]{hsu}) 
 When $A=2$ and $B=1$, we have $u_k=k$ for each nonnegative integer $k$,
and if in addition we assume that $q=p$ is a prime, then  $w_p=p$,   
and hence the congruence \eqref{con83} becomes
  $$
{mp+s\choose np+t}\equiv {m\choose n}{s\choose t}\pmod{p},
  $$    
which is in fact, Lucas' theorem.\hfill$\Box$
  \vspace{1mm}

In 2002 H. Hu \cite[p. 291, Theorem]{hu} proved the following result: 
{\it Let $q$ be a positive integer, and let $m$ and $n$ be even 
nonnegative integers with $n\le m$. Let $s$ and $t$ be nonnegative 
integers such that $t\le s<q$, and let $v_q^{*}$ be the largest divisor 
of $v_q$ relatively prime to $v_0,\ldots,v_{q-1}$. Then
     \begin{equation}\label{con85}
{m/2 \choose n/2}{mq+s\brack nq+t}_u\equiv {m\choose n}{s\brack t}_u
(-B^q)^{\frac{m-n}{2}(nq+t)+\frac{n}{2}(s-t)}\pmod{v_q^{*}}.
  \end{equation}}

Lucas type congruences modulo $p^2$ and $p^3$ ($p$ is a prime $>3$) for 
Lucas $u$-nomial coefficients and Fibonomial  coefficients
are established  in \cite{kw2}, \cite{kw1} and \cite{shi}. Namely,
in 1993 W.A. Kimball and W.A. Webb   \cite{kw2}
(also see \cite[p. 1029]{shi}) proved the following two results: 
{\it Let $p$ be an odd prime and let $m$ and $n$ be nonnegative integers.
Suppose that $\tau$ is the period of the Fibonacci sequence $(F_n)_{n\ge 0}$
modulo $p$, $r$ is the rank of apparition of $p$
$($that is, $r$  is the least index $k$ for which $p$ divides $F_k$$)$,
and $t=\tau/r$ is an integer. In \cite{tf} it is shown
that $t\in\{1,2,4\}$. The number $\varepsilon$ is defined as follows:
$\varepsilon=1$ if $\tau =r$; $\varepsilon=-1$ if $\tau =2r$;
and $\varepsilon^2\equiv -1(\bmod{\, p^2})$ if $\tau =4r$;
in this case $p\equiv 1(\bmod{\, 4})$.
Then
  \begin{equation}\label{con86}
{m\tau \brack n\tau}_{\mathcal F}\equiv {mt \choose nt }
\pmod{p^2}
   \end{equation}
and
   \begin{equation}\label{con87}
{mr\brack nr}_{\mathcal F}\equiv 
\varepsilon^{(m-n)nr} {m\brack n}_{\mathcal F}\pmod{p^2}.
   \end{equation}}

In 1995  Kimball and  Webb \cite[Theorems 1 and 3]{kw1}  proved the 
following results: 
{\it Let $(u_n)_{n\ge 0}$ and $(v_n)_{n\ge 0}$ be the sequences 
defined by \eqref{con76} and \eqref{con77}, respectively, where
$A$ and $B$ are nonzero integers such that $\gcd(A,B)=1$.
Let $p$ be an odd prime, let $\tau$ be the period of 
the sequence $(u_n)_{n\ge 0}$ modulo $p$, and let 
$r$ be the rank of apparition of $p$. Then for 
all   nonnegative integers $m$ and $n$ such that
$n\le m$ there holds
  \begin{equation}\label{con88}
{mr\brack nr}_{u}\equiv 
\left(\frac{v_r}{2}\right)^{(m-n)nr} {m\choose n}\pmod{p^2}
   \end{equation}
and 
 \begin{equation}\label{con87+}
{m\tau \brack n\tau}_{u}\equiv 
\left(1+\frac{1}{2}\tau (m-n)n((-B)^{\tau}-1)\right) {mt\choose nt}\pmod{p^2}.
   \end{equation}}
As a consequence  of the congruence \eqref{con88}, it is proved in 
\cite[Corollary 2]{kw1} that 
   \begin{equation}\label{con88+}
{m\tau\brack n\tau}_{u}\equiv 
\left(1+\tau(m-n)n\left(\left(\frac{v_r}{2}\right)^t-1\right)\right)
 {mt\choose nt}\pmod{p^2}.
   \end{equation}
  Moreover, the congruence \eqref{con87+} immediately implies 
\cite[Corollary 4]{kw1} that {\it if $B=\pm 1$, then}
 \begin{equation}\label{con88++}
{m\tau \brack n\tau}_{u}\equiv {mt\choose nt}\pmod{p^2}.
   \end{equation}

Kimball and Webb \cite[Theorem 5]{kw1} also proved the 
following congruences for the Gaussian $q$-nomial coefficients:
 \begin{equation}\label{con88+++}\begin{split} 
{mr \brack nr}_q &\equiv \left(\frac{q^r+1}{2}\right)^{(m-n)nr}
 {m\choose n}\pmod{p^2}\\
&\equiv \left(1+\frac{1}{2}r(m-n)n(q^r-1)\right){m\choose n}\pmod{p^2},
   \end{split}\end{equation}
{\it where $p$ is a prime, $q$ is any $p$-integral rational number such that 
$q^2-q$ is not divisible by $p$, and $r$ is the rank of apparition of $p$.}
 
In 1998 B. Wilson  \cite{wi} proved the following result: {\it Let
$p$ be a prime such that  $p\not= 2,5$, and let $r$ be the rank  of 
apparition of $p$ with respect to the Fibonacci sequence $(F_n)_{n\ge 0}$.
Then for any nonnegative integers $m,n,s$ and $l$ such that $0\le s,l<r$
  \begin{equation}\label{con89}
{mr\brack nr}_{\mathcal F}\equiv {m\brack n}_{\mathcal F}
F_{r+1}^{(m-n)nr}\pmod{p}
    \end{equation}
and 
  \begin{equation}\label{con90}
{mr +s\brack nr +l}_{\mathcal F}\equiv {m\choose n} {s\brack l}_{\mathcal F}
F_{r+1}^{(nr+l)(m-n)+n(s-l)}\pmod{p}.
    \end{equation}}

In 2007 L.-L. Shi \cite{shi} proved another  
congruence modulo $p^2$ (where $p>3$ is a prime) for 
the Lucas $u$-nomial coefficients. Namely, in \cite[Theorem 2]{shi}   
it is proved the following result:
{\it Let $(u_n)_{n\ge 0}$ be the Lucas sequence defined by \eqref{con76},
where $A$ and $B$ are nonzero integers such that $\gcd (A,B)=1$, and 
$A\not=\pm 1$ or $B\not=1$. Let $p>3$ be a prime not dividing $B$. 
If $r$ is the rank  of apparition of $p$ with respect to 
$(u_n)_{n\ge 0}$, then for 
any nonnegative integers $m,n,s$ and $t$ such that $0\le s,l<r$,  we have  
          \begin{equation}\label{con92}
{mr+s\brack nr+l}_{u}\equiv \left\{\begin{array}{ll}
(-1)^{l-s-1}B^{-{l-s\choose 2}}u_{(m-n)r}u_{l-s}^{-1} & \\
 \times u_{r+1}^{(m-n)(l-1)-n(l-s)}{mr\brack nr}_{u}
\left({l\brack s}_{u}\right)^{-1}  \pmod{p^2} & {\rm if} \quad s<l\\
u_{r+1}^{ml+ns-2nl}\frac{S_{m,s}}{S_{n,l}S_{m-n,s-l}}{mr\brack nr}_{u}
{s\brack l}_{u} \pmod{p^2} & {\rm if} \quad s\ge l,\end{array}\right.
       \end{equation}
where $S_{k,i}=1-(kBu_r)/u_{r+1}\sum_{j=1}^i(u_{j-1}/u_j)$.

If $\Delta:=A^2-4B$ is not divisible by $p$, then ${mr\brack nr}_{u}$
in \eqref{con92} can be replaced by $(v_r/2)^{(m-n)nr}{m\choose n}$.}

In 1995 Kimball and Webb \cite[Theorem]{kw3} and in 2007  L.-L. Shi \cite{shi}
considered the {\it generalized Lucas $u$-nomial coefficients}  and the 
{\it generalized Fibonomial coefficients} defined as follows. 
If $(u_n)_{n\ge 0}$ is the Lucas sequence defined by \eqref{con76} such that 
$A\not=\pm 1$ or $B\not= 1$, and let $(F_n)_{n\ge 0}$ be the Fibonacci 
sequence. For any positive integer $j$ we set 
  $$ 
[n]_u^j=\prod_{k=1}^nu_{kj}\quad {\rm and}\quad 
[n]_{\mathcal F}^j=\prod_{k=1}^nF_{kj},
 $$
for $n=0,1,2,\ldots$, and regard an empty product as value 1.

Then for $n,k=0,1,2,\ldots$ the {\it generalized Lucas $u$-nomial coefficient}
${n\brack k}_u^j$ and the {\it generalized Fibonomial coefficient}  
${n\brack k}_{\mathcal F}^j$ are defined as follows:
    $$
{n\brack k}_u^j=\left\{\begin{array}{ll}
\frac{[n]_u^j}{[k]_u^j[n-k]_u^j} & {\rm if}\quad 0\le k\le n\\
 0  & {\rm otherwise},\end{array}\right.
   $$
     $$
{n\brack k}_{\mathcal F}^j=\left\{\begin{array}{ll}
\frac{[n]_{\mathcal F}^j}{[k]_{\mathcal F}^j[n-k]_u^j} & {\rm if}\quad 0\le k\le n\\
 0  & {\rm otherwise}.\end{array}\right.
   $$
where $(u_{ij}/u_j)_{i\ge 0}$ is also a Lucas sequence.

   In 1995 Kimball and Webb \cite[Theorem]{kw3} extended the congruence 
\eqref{con87} by showing that  
{\it if the rank $r$ of apparition of $p$ is $p+1$ or $p-1$,
then for any prime $p>3$ and any $m\ge n\ge 0$,
    \begin{equation}\label{con91}
{mr\brack nr}_{\mathcal F}\equiv (\mp)^{(m-n)n}{m\brack n}_{\mathcal F}^r
\pmod{p^3}, \quad respectively.
    \end{equation}}

In 2007 Shi \cite{shi} proved  
the congruence modulo $p^3$ (where $p>3$ is a prime) for 
the generalized Lucas $u$-nomial coefficients. 
Namely, in \cite[Theorem 1]{shi}   
it is proved the following result:
{\it Let $A$ and $B$ be nonzero integers such that $\gcd (A,B)=1$, and 
$A\not=\pm 1$ or $B\not=1$. Let $p>3$ be a prime not dividing $B$. 
If the rank $r$ of apparition of $p$ is $p+1$ or $p-1$
$($and hence $r=p-\left(\frac{A^2-4B}{p}\right)$$)$, where 
$\left(\frac{\cdot}{p}\right)$ denotes the Legendre symobol, 
then for any nonnegative integers $m$ and $n$ we have  
     \begin{equation}\label{con93}
{mr\brack nr}_{u}\equiv (-1)^{(m-n)n}B^{(m-n)n{r\choose 2}}{m\brack n}_{u}^r
\pmod{p^3}.
    \end{equation}}

 \vspace{1mm}
{\it Remark} 31. In the case $A=-B=1$ the congruence \eqref{con93} yields 
the congruence \eqref{con91} of  Kimball and Webb \cite{kw3}.\hfill$\Box$
   \vspace{1mm}

In 1965  G. Olive \cite{ol} (also see \cite[Lemma 2.1]{pa2}) proved the 
following result: {\it Suppose that $d$ is a positive integer and
 $a,b,h,l$ are integers such that $0\le b,l\le d-1$. Then 
    \begin{equation}\label{con94}
{ad+b\brack hd+l}_q\equiv {a\choose h}{b\brack l}_q\pmod{\Phi_d(q)},
\end{equation}
where $\Phi_d(q)$ is the $d$th cyclotomic polynomial}.
 
 \vspace{1mm}
{\it Remark} 32. As noticed in \cite[Chapter 5, p. 506]{sc},
the  congruence \eqref{con94} perhaps was known to Gauss and it is
rediscovered  in  1982 by J. D\'esarm\'{e}nien \cite{des} 
and  V. Strehl \cite{str} whose proof uses combinatorial arguments.\hfill$\Box$

{\it Remark} 33. Another different $q$-analogue 
of the congruence \eqref{con94} was established in 1967 by  R.D. Fray \cite{fr2}.\hfill$\Box$
 \vspace{1mm}

{\it Remark} 34. Applying Lucas' theorem, in  2006 S.-P. Eu, S.-C. Liu and Y.-N. 
Yeh \cite{ely} established the congruences of several combinatorial
numbers, including Delannoy  numbers and a class of Ap\'{e}ry-like numbers, 
the numbers of noncrossing connected graphs (Sloane's sequence 
A007297), the numbers of total edges of 
all noncrossing connected graphs on $n$ vertices (Sloane's sequence 
A045741), etc.\hfill$\Box$
  
\section{Some applications of Lucas' theorem}\label{sec6}

Even today, Lucas' theorem is being studied widely, and has both
extended and generalized, particularly in the area of 
divisibility of binomial coefficients. Numerous results 
on divisibility of binomial and multinomial coefficients by primes 
and prime powers and related historical notes are given in 1980  by D. 
Singmaster \cite{sin2}. Furthermore,
Lucas' theorem has numerous applications  in Number Theory, 
Combinatorics, Cryptography and Probability. We also point out that this 
theorem has become ubiquitous in the Theory of cellular automata. 
 
\subsection{Lucas' theorem and the Pascal's triangle}\label{subsec6.1}

Let $a_k(n)$ be the number of integers $0\le m\le n$ such that 
${n\choose m}\not\equiv 0(\bmod{\, k})$, that is, $a_k(n)$
is the number of nonzero entries on row $n$ of Pascal's triangle
modulo $k$. Let $|n|_w$ be the number of occurrences of the word $w$
in $n_sn_{s-1}\cdots n_0$,  where $n=\sum_{i=0}^sn_ik^i$ is the base-$k$ 
representation of $n$.  In 1899 J.W.L. Glaisher \cite[\S 14]{gl} 
initiated the study of counting entries on row $n$ of 
Pascal's triangle modulo $k$ by using Lucas' theorem to determine 
$a_2(n)=2^{|n|_1}$. The proof is simple (cf. \cite[p. 1]{ro}):
In order that ${n\choose m}$ be odd, each term ${n_i\choose m_i}$ 
in the product must be 1, so if $n_i=0$ then $m_i=0$ and if $n_i=1$
then $m_i$ can be either 0 or 1. 
It was the first result on a thorny  path of solution of this difficult 
problem. However, this topic was forgotten for almost a half-century.

In 1947 N.J. Fine \cite{fi} generalized Glaisher's result to an arbitrary 
prime. Fine's result 
follows from Lucas' theorem in the same way: {\it Let $p$ be a prime,
and let $n$ be a nonnegative integer. The number of nonzero entries on row 
$n=\sum_{i=0}^sn_ip^i$ 
of Pascal's triangle modulo $p$ is $($cf. \cite[p. 2]{ro}$)$}
   \begin{equation}\label{con95}
a_p(n)=\prod_{i=0}^{s}(n_i+1).
   \end{equation}
Namely, the formula \eqref{con95} immediately follows from the fact that by
Lucas' theorem, the binomial coefficient ${n\choose m}$ with 
$m=\sum_{i=0}^sm_ip^i$ is not divisible by a prime $p$ if and only if 
$0\le m_i\le n_i$ for all $i=0,1,\ldots,s$.

   \vspace{1mm}
{\it Remark.} 35. If $p=2$, then the formula \eqref{con95} presents
the number of odd entries on row $n=\sum_{i=0}^sn_i2^i$ 
of Pascal's triangle. Notice that the parity of binomial coefficients has 
played an important role in a paper from  1984  of J.P. Jones and Y.V. Matijasevi\v{c}
\cite{jm1} in connection with Hilbert's tenth problem, G\"{o}del's 
undecidability proposition and computational complexity. They base 
their Lemma on the Lucas' theorem given by the congruence 
\eqref{con1} with $p=2$ (cf. \cite[Lemmas 3.9 and 3.10]{jm2}). \hfill$\Box$
   \vspace{1mm}

As noticed in \cite{ro}, one may generalize Glaisher's result in a different
direction, namely to ask for the number $a_{k,r}$ of integers $0\le m\le n$
such that ${n\choose m}\equiv r(\bmod{\, k})$. 
  In 2011 E. Rowland \cite[Section 2, Theorem 1]{ro} generalized Fine's result 
to prime powers, obtaining a formula for the sum
$a_{p^{\alpha}}(n)=\sum_{r=1}^{p^{\alpha}-1,r}(n)$.
Notice that in 1978 E. Hexel and H. Sachs \cite[\S 5]{hs}
determined a formula for $a_{p,r^i}(n)$ in terms of $(p-1)$th roots of unity,
where $r$ is a primitive root modulo $p$.  
For some related results see  also
\cite{as}, \cite{dw3}, \cite{gw},  \cite{gr1} and \cite[Theorem 2]{ro}).

The previous considerations can be genearlized  as follows.
Let $p$ be a prime. For nonnegative integers $n$ and $k$ consider the set 
  $$
A_{n,k}^{(p)}=\{j\in\{0,1,\ldots,n\}: p^k\Vert {n\choose j}\},
  $$
where $p^k\Vert {n\choose j}$ denotes that $p^k\mid {n\choose j}$
and ${n\choose j}\not\equiv 0 (\bmod{\,p^{k+1}})$.
In particular, $A_{n,0}^{(p)}$ is a set of nonzero entries 
on row $n$ of Pascal's triangle modulo $k$. Therefore, 
under the previous notation, for a prime $p$ we have
 $a_p(n)=|A_{n,0}^{(p)}|$
($|S|$ denotes the cardinality of a set $S$), 
Notice that $|A_{n,0}^{(p)}|$ can be evaluated by Fine's formula 
\eqref{con95}. 
In 1967 L. Carlitz \cite{ca} solved a difficult problem for evaluation 
of $|A_{n,1}^{(p)}|$. In   1971  F.T. Howard \cite{ho2}, 
discovered  the formula for $|A_{n,k}^{(2)}|$ for arbitrary $k$.
In 1973 F.T. Howard \cite{ho3}  found a solution for $|A_{n,2}^{(p)}|$.

Further related results are given in \cite{gr}, and in 1997 by 
J.G. Huard, B.K. Spearman and K.S. Williams \cite{hsw}.  
Let $n$ be a nonnegative integer. The $n$th
row of Pascal's triangle consists of the following $n+1$
binomial coefficients:
  $$
{n\choose 0},{n\choose 1},{n\choose 2},\ldots,{n\choose n}.
  $$
We denote by $N_n(t,m)$ the number of those binomial coefficients which 
are congruent to $t$ modulo $m$, where $t$ and $m\ge 1$ are integers  
such that $0\le t\le m-1$. Let  $p$ be a prime, and  
let $n$ be a positive integer with  the  $p$-adic expansion 
$n=\sum_{i=0}^kn_ip^i$. We denote the number of $r$'s occuring among 
$n_0,n_1,\ldots,n_k$ by $l_r$ ($r=0,1,\ldots,p-1$). Set 
$\omega= e^{2\pi i/(p-1)}$ and let $g$ denote a primitive root modulo $p$.
Denote by ${\rm ind}_gt$ the index of the integer $t\not\equiv 0(\bmod{\,p})$ 
with respect to $g$; that is, ${\rm ind}_gt$ is the unique integer $j$ 
such that $t\equiv g^j(\bmod{\, p})$.   In  1978  E. Hexel and H. Sachs
\cite[Theorem 3]{hs} have shown that {\it for $t=1,2,\ldots,p-1$,
     \begin{equation}\label{con102+}
N_n(t,p)=\frac{1}{p-1}\sum_{s=0}^{p-2}\omega^{-s{\rm ind}_gt}
\prod_{r=1}^{p-1}B(r,s)^{l_r},
   \end{equation}
where 
     $$
B(r,s)=\sum_{j=0}^r\omega^{s {\rm ind}_g{r\choose j}}.
     $$} 
By using the formula \eqref{con102+}, in 1997 J.G. Huard, 
B.K. Spearman and K.S. Williams proved the analogous formula for 
$N_n(tp,p^2)$ with $t=1,2,\ldots,p-1$ \cite[Theorem 1.1]{hsw}. 
They proved that {\it for $t=1,2,\ldots,p-1$,
  \begin{equation}\label{con97}
 \begin{split} 
N_n(tp,p^2)= &\frac{1}{p-1}\sum_{i=0}^{p-2}\sum_{j=1}^{p-1}l_{ij}
\sum_{s=0}^{p-2}\omega^{-s({\rm ind}_gt+{\rm ind}_g(i+1)-{\rm ind}_gj)}\\
\times & B(p-2-i,-s)B(j-1,s)\prod_{r=1}^{p-1}B(r,s)^{l_r-\delta(r-i)- 
\delta(r-j)},
   \end{split}\end{equation}
where 
  $$
\delta(x)=\left\{
\begin{array}{ll}
1 & if \,\, x=0\\
0 & if \,\, x\not= 0,
\end{array}\right.
  $$
and $l_{ij}$ denotes the number of occurences of the pair $ij$ in the string 
$n_0n_1\ldots n_k$.}
 
Let $p$ be a prime, and let $k$ be a positive integer. 
 Let  $A(k,p)$ be the matrix with entries 
${i\choose j}_p:={i\choose j}(\bmod{\,p})$, $0\le i<p^k$, $0\le j<p^k$
(actually, ${i\choose j}_p$ is the remainder of the division of 
${i\choose j}$ by $p$). By using the Lucas property of the 
matrix $A(k,p)$ given by \eqref{con54}, in 1994 M. Razpet \cite[p. 378]{raz1} 
{\it proved that the number of all zero entries of the matrix $A(k,p)$
is equal to $p^{2n}-{p+1\choose 2}^k$, and hence,  
the number of all nonzero entries of the matrix $A(k,p)$
is equal to ${p+1\choose 2}^k$}.

Let $p$ be a prime, and let $n$ be  a positive integer.  
For an integer $r$ such that $0\le r\le p-1$, let 
$b_r(n)$ be the number of binomial coefficients 
${i\choose j}$ with  $0\le j\le i<n$ such that
${i\choose j}\equiv r(\bmod{\, p})$.
In 1957 J.B. Roberts \cite{r}  established systems
of simultaneous linear difference equations with 
constant coefficients whose solutions would yield the quantities 
$b_r(n)$ explicitly.  Namely, {\it if $0\le c\le p-1$, $1\le t\le p^k$, $k>0$,
and if $\bar{q}$ is the reciprocal  of $q\in\{1,3,\ldots,p-1\}$
modulo $p$   $($i.e., $q\bar{q}\equiv 1(\bmod{\,p})$$)$, then 
by \cite[Theorem 1]{r},
   \begin{equation}\label{con98}
b_r(cp^k+t)=b_r(cp^k)+\sum_{q=1}^{p-1}(b_{r\bar{q}}(c+1)-b_{r\bar{q}}(c))
b_q(t).
  \end{equation}}
{\it Furthermore, if $b(n)=\sum_{r=1}^{p-1}b_r(n)$ and 
$n=\sum_{i=0}^kn_ip^i$ with $0\le n_i\le p-1$ for all $i=0,1,\ldots,k$, then 
by \cite[Corollary 4]{r},
    \begin{equation}\label{con99}
b(n)=\frac{1}{2}\sum_{i=0}^kn_i((n_i+1)\cdots(n_k+1))
\left(\frac{1}{2}p(p+1)\right)^i.
   \end{equation}} 

By using Lucas' theorem, in 1992 R. Garfield and H.S. Wilf \cite[Theorem]{gw}
proved the following result: 
{\it Let $p$ be a prime, let $a$ be a primitive root modulo $p$,
and let $n$ be a nonnegative integer with the $p$-adic expansion
$n=\sum_{i=0}^sn_ip^i$. Denote by $l_j(n)$ 
the number of $j$'s occuring among  
$n_0,n_1,\ldots,n_s$ $($$j=0,1,\ldots,p-1$$)$. Further, for each 
$i\in\{0,1,\ldots ,p-2\}$ 
let $r_i(n)$ be the number of integers $k$ with $0\le k\le n$, for which
${n\choose k}\equiv a^i(\bmod\,p)$, and let 
$R_n(X)=\sum_{i=0}^{p-2}r_i(n)X^i$ be their generating function. 
Then 
     \begin{equation}\label{con100}
R_n(X)\equiv \prod_{j=1}^{p-1}R_j(X)^{l_j(n)}\pmod{X^{p-1}-1}.
   \end{equation}} 

In 1990 R. Bollinger and C. Burchard \cite{bb} considered the 
extended pascal's triangles which arise, by analogy with the ordinary 
Pascal's triangle as the (left-justified) arrays of the coefficients 
in the expansion $(1+x+x^2+\cdots +x^{k-1})^n$. That is, 
the array $T_k$ has in row $n$, column $m$, 
the  number $C_k(n,m)$ defined for $k,n,m\ge 0$ by the expansion
   $$
(1+x+x^2+\cdots +x^{k-1})^n=\sum_{m=0}^{(k-1)n}C_k(n,m)x^m,
 $$
It is nociced in \cite[the property d) on page 199]{bb}
that 
  $$
C_k(n,m)=\sum_{j}(-1)^j{n\choose j}{n-1+m-kj\choose n-1},
  $$
and hence, $C_2(n,m)={n\choose m}$. Accordingly, $T_2$ is the 
Pascal's triangle.

R. Bollinger and C. Burchard \cite[Theorem 1]{bb} applied Lucas' theorem to the  Pascal's triangle,
proving that {\it if $p$ is a prime, and if
$n=n_0+n_1p+\cdots +n_sp^s$ and
$m=m_0+m_1p+\cdots +m_sp^l$ are the $p$-adic expansions of $n$ and $m$,
then
 \begin{equation}\label{con101}
C_k(n,m)\equiv \sum_{r_0,\ldots,r_s}\prod_{i=0}^sC_{k}(n_i,r_i)\pmod{p}, 
 \end{equation}
where the sum is taken over all 
$(s+1)$-tuples $(r_0,r_1,\ldots,r_s)$ such that  
$i)$ $m=r_0+r_1p+\cdots +r_sp^s$ and $ii)$ $0\le r_i\le (k-1)n_i$ for each 
$i=0,1,\ldots,s$; if $m$ is not representable in this form, then certainly
$C_k(n,m)\equiv 0\,(\bmod{\,p})$.}

\subsection{Another applications of Lucas's theorem}\label{subsec6.2}

By using Kummer's theorem and Lucas' theorem, in 2007 K. Dilcher 
\cite[Theorem 2]{di} derived an alternating sum analog to a special 
case to an 1876 congruence of Hermite \cite{he} (also see 
\cite[Chapter IX, p. 271]{d}) as follows. {\it Let $p$ be an odd prime and 
let $q$ be a positive integer. Then}
     \begin{equation}\label{con102}   
\sum_{j=0}^{\lfloor q/2 \rfloor}{q(p-1)\choose 2j(p-1)}\equiv \left\{
    \begin{array}{ll}
1 \pmod{p} & if\,\, q \,\, is \,\, odd;\\
2 \pmod{p} & if\,\, q \,\, is \,\, even\,\, and \,\, q\not\equiv 0 
(\bmod\,p+1);\\
\frac{3}{2} \pmod{p} & if\,\, p+1\mid q.
  \end{array}\right.
    \end{equation}

By using Lucas' theorem, in 2009 the author of this article proved the 
following result \cite[Theorem]{me}. 
{\it If $d,q>1$ are integers such that 
  \begin{equation}\label{con103}   
{nd\choose md}\equiv {n\choose m} \pmod{q}
       \end{equation}
for every pair of integers $n\ge m\ge 0$,
then $d$ and $q$ are powers of the same prime $p$}.

   \vspace{1mm}
{\it Remark} 36.  Observe that the above result 
may be considered as a partial converse theorem of 
the congruence \eqref{con5}  of Subsection \ref{subsec2.1}.\hfill$\Box$
   \vspace{1mm}

In 2010 M.P. Saikia and J. Vogrinc 
\cite[Theorem 2.1]{svo} (see also 
\cite[Theorem 1.2 and its proof]{ls})  proved that {\it a positive
integer $p>1$ is a prime if and only if 
       \begin{equation}\label{con104}   
{n\choose p}\equiv \left\lfloor\frac{n}{p}\right\rfloor\pmod{p}
       \end{equation} 
for every nonnegative integer $n$.}
  
By using Lucas' theorem, in 2013 the author of this article 
\cite[Theorem 1.1]{me7} generalized Babbage's criterion 
for primality given in 1819 by Babbage \cite{bab}
(also see \cite[Section 4]{gr}). Lucas' theorem is also applied 
in a recent author's note \cite[Theorem 1]{me9} in order 
to prove the following  result:
{\it If $n>1$ and $q>1$ are integers such that 
  $$
{n-1\choose k}\equiv (-1)^k \pmod{q}
   $$
for every integer $k\in\{0,1,\ldots, n-1\}$,
then $q$ is a prime and $n$ is a power of $q$.}

\vspace{1mm}

 {\it Definition} (see, e.g., \cite{ags}). Let $p$ be a prime. We say that the 
sequence of rational numbers $(a_n)_{n\ge 0}$ $(a_n)_{n\ge 0}$ has the 
$p$-{\it Lucas property}  
(or that the sequence $(a_n)_{n\ge 0}$ is  $p$-{\it Lucas}) 
if the denominators of all the $a_n$'s 
are not divisible by $p$, and if for all $n\ge 0$ and for all $j\in \{0,1,\ldots,p-1\}$ it holds
    \begin{equation}\label{con105}
\qquad \qquad \qquad  a_{pn+j}\equiv a_na_j\pmod{p}.\qquad  \qquad\qquad  
\qquad  
\hfill\Box
    \end{equation}
Clearly, the sequence of rational numbers 
$(a_n)_{n\ge 0}$ has the $p$-Lucas property if and only if 
     \begin{equation}\label{con106}
a_n\equiv \prod_{i=0}^sa_{n_i}\pmod{p},
    \end{equation}
for every positive integer $n$ with the $p$-adic expansion
 $n=n_0+n_1p+\cdots +n_sp^s$ such that  $0\le n_i\le p-1$ for all 
$i=0,1,\ldots,s$. Furthermore, the integer sequence 
$(a_n)_{n\ge 0}$  has the Lucas property if and only if 
$(a_n)_{n\ge 0}$  has the $p$-Lucas property for every prime $p$. 

In what follows, we will consider sequences 
$(a_n)_{n\ge 0}$ having the $p$-Lucas property for infinitely many primes $p$. 
As noticed in \cite[Remarks 1]{ags}, such 
a sequence is either 0 or it satisfies $a_0=1$.\hfill$\Box$

For a positive integer $t$ consider {\it the formal power series}
   $$
\sum_{n=0}^{\infty}{2n\choose n}^tX^n.
  $$ 
It is known that the above formal power series is transendental over 
$\Bbb Q(X)$ when $t\ge 2$. This is due in 1980 to Stanley \cite{st2},
and independently in 1987 to Flajolet \cite{fl} and in 1989 to  
C.F. Woodcock and H. Sharif \cite{ws}. While Stanley and Flajolet used 
analytic methods and studied the asymptotics of the coefficients 
of this series, Woodcock and  Sharif gave a purely algebraic proof.
Their basic idea is to reduce this series modulo a prime $p$, and to
use the $p$-Lucas property for central binomial coefficients: if 
$n=\sum_{i=0}^sn_i$  is the base $p$ expansion of a positive integer $n$,
then (\cite{m2}; cf. \eqref{con57} of Subsection \ref{subsec4.1}) 
       \begin{equation}\label{con107}
{2n\choose n}\equiv \prod_{i=0}^s{2n_i\choose n_i}\pmod{p}.
      \end{equation}
Namely, a proof of  Woodcock and Sharif \cite{ws}  is based on 
the following congruence which follows from Lucas' theorem:
    $$
F_t^{p-1}(X)\equiv \left(\sum_{i=0}^{(p-1)/2}{2i\choose i}X^i\right)^{-1}
\pmod{p}.  
     $$
In 1998 J.-P. Allouche, D. Gouyou-Beauchamps and G. Skordev 
\cite{ags}
generalized the method of Woodcock and  Sharif to characterize 
all formal power series that have the $p$-Lucas property for 
``many'' primes $p$, and that are furthermore algebraic over 
$\Bbb Q(X)$. Namely, they proved the following result \cite[Theorem 1]{ags}:
 {\it Let $s$ be an integer $\ge 2$. Define $s'=s$ if $s$ is even, and 
$s'=2s$ if $s$ is odd. Let $F(X)=\sum_{n=0}^{\infty}a_nX^n$ be a nonzero
formal power series with coefficients in $\Bbb Q$. 
Then the following conditions are equivalent:
  \begin{itemize}
\item[(i)] The sequence $(a_n)_{n\ge 0}$ has the $p$-Lucas property for 
all large primes $p$ such that $p\equiv 1(\bmod{\,s})$, and the 
formal power series $F(X)$
is algebraic over $\Bbb Q(X)$.
\item[(ii)] There exists a polynomial $P(X)$ in $\Bbb Q[X]$ of degree 
at most $s'$, with $P(0)=1$, such that 
$F(X)=(P(X))^{-1/s'}$.
 \end{itemize}

If $s$ is odd, and if the number $s'$ is replaced by $s$ in the statement 
$(ii)$, we still have $(ii)$ implies $(i)$, but the converse is not 
necessarily true.}    

Furthermore, when the number $s$ is equal to 2,
in 1999 Allouche \cite[Theorem 6.4]{al1} proved the following 
result (cf. \cite[Theorem 2]{ags}):
{\it Let $(a_n)_{n\ge 0}$ be a nonzero sequence of rational numbers. Then
the following assertions are equivalent.
  \begin{itemize}
\item[(i)] The sequence $(a_n)_{n\ge 0}$ has the $p$-Lucas property for 
all large primes $p$, and the series $F(X)=\sum_{n=0}^{\infty}a_nX^n$
is algebraic over $\Bbb Q(X)$.

\item[(ii)] For all large primes $p$ the sequence $(a_n)_{n\ge 0}$ has the
$p$-Lucas property, and the degree $d_p$ of the series 
$\sum_{n=0}^{\infty}(a_n(\bmod{\, p}))X^n$ $($that is necessarily algebraic
over $\Bbb F_p(X)$ from the $p$-Lucas property$)$ is bounded 
independently of $p$. 

\item[(iii)] There exists a polynomial $P(X)$ in $\Bbb Q[X]$ of degree 
at most $2$, with $P(0)=1$, such that 
$F(X)=\sum_{n=0}^{\infty}a_nX^n=(P(X))^{-1/2}$.
 \end{itemize}}    

   \vspace{1mm}
{\it Remark} 37. In 2013 \'{E}. Delaygue \cite[Subsection 1.2]{de}
considered the notion of $p$-Lucas property to a 
$\Bbb Z_p$-{\it valued family} 
$A=\left(A(\mathbf{n})\right)_{\mathbf{n}\in \Bbb N^d}$, 
where $p$ is a prime, $\Bbb Z_p$ is the ring of $p$-adic integers
and  $d$ is a positive integer.  We say that $A$  satisfies the
$p$-Lucas property if and only if, for all 
$\mathrm{v}\in \{0,1,\ldots,p-1\}^d$ and all $\mathrm{n}\in \Bbb N^d$,
we have 
   $$
A(\mathrm{v}+\mathrm{n}p)\equiv A(\mathrm{v})A(\mathrm{n})\pmod{p\Bbb Z_p}. 
  $$
Delaygue \cite[Theorem 3]{de} established an effective criterion 
for a sequence of factorial ratios to satisfy the $p$-Lucas property for 
almost all primes $p$.\hfill$\Box$

\vspace{4mm}

\section*{Appendix}

\centerline{\bf{List of references and related  Lucas type congruences}} 

\centerline{ \bf{from this article (arranged by year of publication)}}

\vspace{2mm}

{\small

\cite[1819]{bab} Charles Babbage    (also 
see \cite[Introduction]{gr} or \cite[page 271]{d}) - \eqref{con17}, p. 10.

\cite[1862]{w} Joseph Wolstenholme - \eqref{con18}, p. 10.  

\cite[1869]{an}
 H. Anton (also see \cite[p.  271]{d}) - \eqref{con2}, p. 5)

\cite[1878]{l}, \cite[1878; Section XXI, pp. 229--230]{l2},  \'{E}. Lucas 
(Lucas' theorem) - the congruences \eqref{con1} and \eqref{con4}, p. 5.

 \cite[1900; p. 323]{gl2} 
J.W.L. Glaisher (also see \cite[the congruence 7.1.5]{no}  
and  \cite[Section 6]{me5}) - \eqref{con19}, p. 10.

\cite[1949]{bs} W. Ljunggren (also see \cite[Theorem 4]{ba}, \cite{gr}, 
\cite[Problem 1.6 (d)]{st} and \cite{sio}) - \eqref{con20}, p. 11.

  \cite[1952]{bs} 
E. Jacobsthal (also see \cite{gr}) -   \eqref{con21}, p. 11.

\cite[1955]{ca2} L. Carlitz   -   \eqref{con69}, p. 25.

 \cite[1956]{cl} L.E. Clarke
and  \cite[1957]{pi}, P.A. Piza - \eqref{con7}, p. 8.

\cite[1965]{ol}  G. Olive (also see \cite[Chapter 5, p. 506]{sc}, 
\cite{des}, \cite{str} \cite[Lemma 2.1]{pa2}) - \eqref{con94}, p. 36.

\cite[1982; Theorem 1]{ge}  I. Gessel   - \eqref{con50}, p. 19.

\cite[1987; Theorem 1]{ho4} F.T. Howard - \eqref{con68}, p. 25.

\cite[1987; p. 306, Corollary and Theorem 2]{ho4} 
F.T. Howard - \eqref{con70} and \eqref{con71}, p. 25.

 \cite[1988; Theorem 2]{ma}, R.A. Macleod - \eqref{con40}, p. 16.

\cite[1988; p. 61, Theorem]{sv} M. Sved - \eqref{con63}, p. 23. 

\cite[1988; p. 60, Theorem]{sv} M. Sved - \eqref{con79}, p. 31. 

 \cite[1990; Theorem 3]{ba} D.F. Bailey - \eqref{con38}, p. 12.

 \cite[1990; Theorem 5]{ba} D.F. Bailey - 
\eqref{con38+}, p. 13.

 \cite[1990; Theorem 1]{bb} R. Bollinger and C. Burchard - 
\eqref{con101}, p. 40.

 \cite[1991; Theorem 4]{ba2} D.F. Bailey - 
\eqref{con23}, p. 11.

 \cite[1991; Theorem 5]{ba2} D.F. Bailey - 
\eqref{con24}, p. 12.

\cite[1991; Theorem 3]{dw1} K. Davis and W. Webb  
(also see  \cite[p. 88, Theorem 5.1.2]{lo} and 
\cite[p. 34, congruence (2.2)]{bh}) - \eqref{con32}, p. 14.

 \cite[1992; Theorem 5]{ba3} D.F. Bailey (also see 
\cite[Corollary 1.2]{me8}) - \eqref{con25} and \eqref{con26}, p. 12.

\cite[1992; Proposition 2]{gr1} A. Granville - \eqref{con29}, p. 13.

\cite[1992]{m2} R.J. McIntosh  -  \eqref{con107}, p. 42.

 \cite[1992; Theorem 2]{we3} 
(also see \cite{we4} and  \cite[Section 7]{hol1}) D.L. Wells  
- \eqref{con82+} and \eqref{con82++}, p. 32.

\cite[1993; Theorem 3]{dw2} K. Davis and W. Webb 
(also see \cite[Proposition 2]{gr1}) - \eqref{con35}, p. 15. 

\cite[1993; Corollary 1]{dw2} K. Davis and W. Webb - \eqref{con36}, p. 15. 

\cite[1993; Corollary 1]{dw2} 
 K. Davis and W. Webb (also see \cite[Theorem]{me4}) - \eqref{con37}, p. 15. 

   \cite[1993]{kw2} W.A. Kimball and W.A. Webb (also see \cite[p. 1029]{shi})
- \eqref{con86} and \eqref{con87}, p. 33. 

 \cite[1993; Proposition 2.1]{prw} R. Peele, A.J. Radcliffe and H.S. Wilf 
 - \eqref{con67}, p. 24.

\cite[1994; p. 60]{hol2}   J.M. Holte (also  see  \cite[1994; p. 227]{hol1}) 
- \eqref{con52} and \eqref{con53}, p. 9.

\cite[1994; p. 378]{raz1} M. Razpet - \eqref{con54}, p. 20.

 \cite[1994; Theorem 2]{we1} D.L. Wells  - \eqref{con78}, p. 30.

\cite[1995; Section 6, the congruence (24)]{gr} 
A. Granville - \eqref{con14}, p. 9.

\cite[1995; Theorem 1]{gr} A. Granville - \eqref{con28}, p. 13.

 \cite[1995; Theorems 1 and 3]{kw1} W.A. Kimball and W.A. Webb 
- \eqref{con88} and \eqref{con87+}, p. 34.

 \cite[1995; Corollaries  2 and 4]{kw1} W.A. Kimball and W.A. Webb 
- \eqref{con88+} and \eqref{con88++}, p. 34.

 \cite[1995; Theorem 5]{kw1} W.A. Kimball and W.A. Webb 
- \eqref{con88+++}, p. 34.

 \cite[1995; Theorem]{kw3} W.A. Kimball and W.A. Webb 
- \eqref{con91}, p. 35.

\cite[1998; Lemma 4]{cal} N.J. Calkin  - \eqref{con55+}, p. 20.

 \cite[1998]{wi} B. Wilson - \eqref{con89} and \eqref{con90}, p. 34.

 \cite[1999; Proposition 7.1]{al1}  J.-P. Allouche - \eqref{con49}, p. 18.

\cite[2000; Theorem 1]{hol1} J.M. Holte - \eqref{con80}, p. 31.

\cite[2000; Theorem 3]{hol1} J.M. Holte - \eqref{con81}
and \eqref{con82}, p. 32.

 \cite[2000; Proposition 3.1]{pe}  R. S\'{a}nchez-Peregrino - 
\eqref{con64} and \eqref{con65}, pp. 23--24.

 \cite[2000; Proposition 4.1]{pe}  R. S\'{a}nchez-Peregrino - 
\eqref{con66}, p. 24.

 \cite[2001; Theorem 2.2]{bc} J. Boulanger and J.-L. Chabert
- \eqref{con75}, p. 27. 

\cite[2001; Theorem]{hsu} H. Hu and Z.-W. Sun - \eqref{con83}
and \eqref{con84}, pp. 32--33. 

 \cite[2002; p. 291, Theorem]{hu} H. Hu  - \eqref{con85}, p. 33.

 \cite[2002; Theorem 1]{raz2}  M. Razpet - \eqref{con61}, p. 22.

\cite[2003; Theorem 5]{bk} D. Berend and N. Kriger  - \eqref{con55}, p. 20. 

 \cite[2004; Theorem 1]{pa1}  H. Pan - \eqref{con62}, p. 23.  

\cite[2005; Theorem 3]{ev}  T.J. Evans - \eqref{con15+}, p. 9. 

  \cite[2006; Theorem 4.7]{ds} E. Deutsch and B.E.  Sagan 
- \eqref{con56}, p. 21. 

  \cite[2006; Theorem 4.4]{ds} E. Deutsch and B.E.  Sagan 
- \eqref{con57}, p. 21. 

\cite[2007; Theorem 2]{shi} L.-L. Shi - \eqref{con92}, p. 35. 

\cite[2007; Theorem 1]{shi} L.-L. Shi - \eqref{con93}, p. 36. 

\cite[2007; Theorem 1.7]{sd}  Z.-W. Sun and D.M. Davis 
- \eqref{con74}, p. 26.  

 \cite[2008; the congruence (1) of Corollary on 
page 490]{ht}  C. Helou and G. Terjanian (also see \cite[Section 11.6, 
Corollary 11.6.22, p. 381]{co}) -  \eqref{con22}, p. 11. 

\cite[2008; (ii) of Lemma 2]{jls} Y. Jin, Z-J. Lu and A.L. Schmidt 
- \eqref{con51}, p. 19.

\cite[2008; Theorem 1.1]{sw} Z.-W. Sun and D. Wan - \eqref{con73}, p. 26. 

\cite[2009; Theorem 2.2]{chd1}  M. Chamberland and K. Dilcher 
- \eqref{con58}, p. 22.

\cite[2009; Corollaries  2.1 and 2.2]{chd1}  M. Chamberland and K. Dilcher 
- \eqref{con59} and \eqref{con60}, p. 22.

\cite[2011; the congruences 7.3.1--7.33]{no} A. Nowicki  - \eqref{con6}, p. 7.

\cite[2012; Theorem 1.1]{me8} R. Me\v{s}trovi\'{c}
- \eqref{con27}, p. 12. 

\cite[2012; Theorem]{me4} R. Me\v{s}trovi\'{c}
(also see \cite[Corollary 1]{dw2}) - \eqref{con37}, p. 15. 

\cite[2012; Proposition]{me4} R. Me\v{s}trovi\'{c}
 - \eqref{con45} and \eqref{con46}, p. 17. 

\cite[2014; Section 5, Theorem 5.3]{ro2} E. Rowland and R. Yassawi -
\eqref{con41}, p. 16. 
 \end{document}